 \documentclass[3p]{elsarticle}
\usepackage{graphicx}
\usepackage{amsmath,amssymb,mathtools,amsfonts}
\usepackage{hyperref}
\usepackage{multirow}
\usepackage{color}
\usepackage{psfrag}
\usepackage{empheq}
\usepackage{tabularx}
\usepackage{bm} 
\usepackage{url}
\usepackage{float}
\usepackage{comment}
\usepackage{placeins}

\usepackage{tikz}
\usepackage{pgfplots}
\usetikzlibrary{arrows.meta}
\newlength\fwidth
\newlength\fheight
\usetikzlibrary{arrows}

\usepackage{algorithm}
\usepackage{algpseudocode}
\usepackage{comment}
\usepackage{caption, subcaption}

\newcommand{\R}{{\mathbb{R}}}
\newcommand{\dd}{{\mathrm{d}}}

\newcommand{\Nstg}{{n_\mathrm{s}}} 
\newcommand{\NFE}{N_\mathrm{FE}} 

\newcommand{\fesd}{{\mathrm{fesd}}}
\newcommand{\Tctrl}{{T}}
\newcommand{\Nctrl}{N}

\newcommand{\zn}{\lambda_{\mathrm{n}}}
\newcommand{\zt}{\lambda_{\mathrm{t}}}

\newcommand{\zttilde}{\tilde{\lambda}_{\mathrm{t}}}
\newcommand{\zttildei}{\tilde{\lambda}_{\mathrm{t}}^{i}}
\newcommand{\zni}{\lambda_{\mathrm{n}}^{i}}
\newcommand{\zti}{\lambda_{\mathrm{t}}^{i}}
\newcommand{\vti}{v_{\mathrm{t}}^{i}}

\newcommand{\Jni}{J_{\mathrm{n}}^{i}}
\newcommand{\Jti}{J_{\mathrm{t}}^{i}}
\newcommand{\fc}{f_{c}}
\newcommand{\fci}{f_{c}^{i}}
\newcommand{\Dti}{D^{i}}
\newcommand{\zdi}{\lambda_{d}^{i}}

\newcommand{\Zn}{\Lambda_{\mathrm{n}}}
\newcommand{\Zt}{\Lambda_{\mathrm{t}}}
\newcommand{\Zni}{\Lambda_{\mathrm{n}}^{i}}
\newcommand{\Zti}{\Lambda_{\mathrm{t}}^{i}}

\newcommand{\Zdin}{\Lambda_{d,n}^{i}}

\newcommand{\Slakcposi}{\Upsilon^{i}_{+}}
\newcommand{\Slakcnegi}{\Upsilon^{i}_{-}}

\newcommand{\Slakcposin}{\Upsilon^{i}_{+,n}}
\newcommand{\Slakcnegin}{\Upsilon^{i}_{-,n}}
\newcommand{\Slakcposn}{\Upsilon_{+,n}}
\newcommand{\Slakcnegn}{\Upsilon_{-,n}}

\newcommand{\Znn}{\Lambda_{\mathrm{n},n}}
\newcommand{\Ztn}{\Lambda_{\mathrm{t},n}}
\newcommand{\Znin}{\Lambda_{\mathrm{n},n}^{i}}
\newcommand{\Ztin}{\Lambda_{\mathrm{t},n}^{i}}

\newcommand{\ts}{t_{\mathrm{s}}}
\newcommand{\tsp}{t_{\mathrm{s}}^+}
\newcommand{\tsn}{t_{\mathrm{s}}^-}

\newcommand{\ncontacts}{n_{\mathrm{c}}}

\newcommand{\nt}{n_{\mathrm{t}}}

\newcommand{\CoR}{\epsilon_{\mathrm{r}}}
\newcommand{\CoRi}{\epsilon_{\mathrm{r}}^{i}}

\newcommand{\sign}{{\mathrm{sign}}}

\newcommand{\FESDj}{FESD-J}

\newtheorem{theorem}{Theorem}

\newtheorem{remark}[theorem]{Remark}

\newtheorem{lemma}[theorem]{Lemma}

\DeclareMathSymbol{\shortminus}{\mathbin}{AMSa}{"39}

\usepackage[shortlabels]{enumitem} 

\begin{document}

\begin{frontmatter}


\fntext[label1]{Corresponding author.} 	
 \fntext[label2]{This research was supported by th DFG via Research Unit FOR 2401 and project 424107692, and by the EU via ELO-X 953348. 
  	\textit{Email address}: \texttt{\{armin.nurkanovic,jonathan.frey,moritz.diehl\} @imtek.uni-freiburg.de, anton.pozharskiy@ merkur.uni-freiburg.de}.}

\title{FESD-J: Finite Elements with Switch Detection for Numerical Optimal Control of Rigid Bodies with Impacts and Coulomb Friction}

\author[FreiburgImtek,label1]{Armin Nurkanovi\'c}
\author[FreiburgImtek,FreiburgMath]{Jonathan Frey}
\author[FreiburgImtek]{Anton Pozharskiy}
\author[FreiburgImtek,FreiburgMath]{Moritz Diehl}
%
\affiliation[FreiburgImtek]{organization={Department of Microsystems Engineering (IMTEK), University of Freiburg}, country={Germany}}
\affiliation[FreiburgMath]{organization={Department of Mathematics, University of Freiburg}, country={Germany}}
%
\begin{abstract}
The Finite Elements with Switch Detection (FESD) is a high-accuracy method for the numerical simulation and solution of optimal control problems subject to discontinuous ODEs. 
In this article, we extend the FESD method \cite{Nurkanovic2022} to the dynamic equations of multiple rigid bodies that exhibit state jumps due to impacts and Coulomb friction. 
This new method is referred to as FESD with Jumps (FESD-J).
Starting from the standard Runge-Kutta equations, we let the integration step sizes be degrees of freedom. 
Additional constraints are introduced to ensure exact switch detection and to remove spurious degrees of freedom if no switches occur. 
Moreover, at the boundaries of each finite element, we impose the impact equations in their complementarity form, at both the position and velocity level. 
They compute the normal and tangential impulses in case of contact making. 
Otherwise, they are reduced to the continuity conditions for the velocities.
FESD-J treats multiple contacts, where each contact can have a different coefficient of restitution and friction. 
All methods introduced in this paper are implemented in the open-source software package NOSNOC~\cite{Nurkanovic2022b}.
We illustrate the use of FESD-J in both simulation and optimal control examples.
\end{abstract}

\end{frontmatter}
\section{Introduction \& related work}\label{sec:introduction}
When perfectly rigid bodies get into contact the normal velocities must immediately jump to a nonnegative value to avoid interpenetration.
Moreover, Coulomb friction models introduce at the same time jumps in the tangential directions. 
Another source of nonsmoothness are frictional slip-stick transitions.
The rigidity and Coulomb friction modeling assumptions simplify the description of the microscopic event of contact and reduce its macroscopic parametrization to only two parameters: the coefficient of restitution and friction.
However, the nonsmoothness complicates the numerical treatment of such systems in simulation and optimal control.

In this paper, we propose an equation-based event-driven method for nonsmooth rigid bodies models, which can be seamlessly used in direct transcription methods of Optimal Control Problems (OCPs).
The method does not need a zero-finding or mode selection procedure, typically used in event-driven methods.
The proposed method has the same accuracy that the underlying Runge-Kutta methods have for smooth ODEs and DAEs.

For simulating such systems the usual choice are (semi-)implicit Euler time-stepping methods that require solving a Linear Complementarity Problem (LCP) at each time step.
These methods are known for their stability, ease of use, and ability to handle a large number of contacts~\cite{Acary2008,Jean1992,Jean1999,Moreau1988,Moreau1999,Anitescu1997,Stewart2000}.
However, they are limited to first-order accuracy, even in the absence of contacts.
Moreover, they approximate the nonlinear friction cone with a polyhedral cone.
One of the first numerical methods in this class is the Jean-Moreau time-stepping scheme~\cite{Jean1992,Jean1999,Moreau1988,Moreau1999}.
It treats the unilateral constraints on velocity level and uses Newton's impact law.
The Schatzman–Paoli~\cite{Paoli2002,Paoli2002a} scheme deals directly with the inequality constraints at the position level.
Similar widely used methods are the Anitescu-Potra method (at velocity level)~\cite{Anitescu1997} and the Stewart-Trinkle method (at position level)~\cite{Stewart1996a}.
A general convergence result of time-stepping methods was provided by Stewart~\cite{Stewart1998}.
Several modifications have been made to these methods, such as those to handle stiff systems~\cite{Anitescu2004} or to improve computational efficiency by solving convex quadratic programs instead of LCPs~\cite{Anitescu2006}.
A comprehensive review of time-stepping methods can be found in~\cite{Studer2008}.
Event-driven methods are less commonly used in the simulation of rigid bodies, cf.~\cite{Acary2008} for an overview.
Moreover, they are not practical for discretizing optimal control problems due to their external switch-detection procedure, which cannot easily be incorporated into a direct transcription method.

In direct methods, a continuous time Optimal Control Problem (OCP) is first discretized and one obtains a finite-dimensional Nonlinear Program (NLP).
The ultimate goal of the methods presented in this paperis to solve optimal control problems subject to nonsmooth rigid body dynamics via direct methods.
In the past two decades, so-called \textit{contact-implicit trajectory optimization} received a lot of attention in the robotics community.
Trajectory optimization problems subject to such dynamics result in nonsmooth OCPs, where the contact sequence of the rigid bodies is discovered fully implicitly. 
Usually, a direct approach is taken to solve these OCPs. 
One of the first works to do so is by Posa and Tedrake~\cite{Posa2014}. 
They use an implicit Euler discretization~\cite{Stewart1996a,Anitescu1997} for the discretization and solve the MPCC with an SQP method. 
Similar approaches with different (possibly smoothed) contact models, adapted time-stepping methods, and MPCC solution strategies are presented in~\cite{Carius2018,Doshi2020,Howell2022a,Manchester2019,Tassa2012}, to name a few.
The recent survey \cite{Wensing2022} provides a broad overview of the use of optimization in legged locomotion.

Numerical sensitivities are crucial for the success of direct optimal control, where the NLPs are solved with Newton-type methods.
However, most time-stepping methods and physics engines based on them provide in general wrong sensitivities~\cite{Zhong2022}.
This is no surprise since the same issues arise in the time-discretization of discontinuous ODEs~\cite{Nurkanovic2022,Stewart2010}.
Besides the low accuracy, this can lead the solver to get stuck at artificial local minima arbitrarily close to the initialization point.
However, the approaches sometimes lead to practically satisfying results, since optimization with wrong derivatives still can yield feasible solutions~\cite{Bock2007}.
Moreover, during the solution process, the discretized OCP is smoothed explicitly or implicitly, which improves the convergence in the early iterates~\cite{Stewart2010,Nurkanovic2020}.

In the case of ODEs with a discontinuous r.h.s., these difficulties can be overcome with the Finite Elements with Switch Detection (FESD) method~\cite{Nurkanovic2022b,Nurkanovic2023,Nurkanovic2022}.
In this method, inspired by Baumrucker and Biegler~\cite{Baumrucker2009}, the integration step sizes $h_k$ are degrees of freedom.
Additional \textit{cross complementarity} conditions are introduced that ensure exact switch detection.
If no switches occur there are spurious degrees of freedom are removed by so-called \textit{step equilibration} conditions.
However, FESD cannot be applied directly to rigid body systems with state jumps.
One solution is to use the \textit{time-freezing} reformulation \cite{Nurkanovic2021,Nurkanovic2021a,Halm2021}, which transforms rigid bodies with friction and impact models into ODEs with a discontinuous r.h.s., with a continuous solution on a different time domain.
The FESD method can be applied to systems obtained by the time-freezing reformulation and recover the exact solution and thereby overcome the issues of low accuracy and wrong sensitivities.
In this paper, we propose an extension to the FESD method, which we call \textit{\FESDj}, that can be directly applied to the nonsmooth systems with state jumps without the use of the time-freezing reformulation.

The works most similar to the present paper are \cite{Patel2019,Shield2022}.
In \cite{Patel2019}, the authors consider an orthogonal collocation discretization for rigid bodies with inelastic impacts and friction.
Inspired by \cite{Baumrucker2009}, they also let the step sizes be degrees of freedom and introduce cross complementarity-type constraints, which allow switching only at the finite element boundaries.
Continuity conditions are imposed on the velocity variables for two adjacent finite elements, and the state jumps are essentially treated within a single finite element that approximates the entire impacting impulse.
The error introduced by this approximation shrinks with element size.
Therefore, the length of this element must be sufficiently close to zero to recover high discretization accuracies.
However, there is no guarantee that the optimizer will always choose an appropriate value.
In \cite{Shield2022}, the same group improves this approach by imposing impulse equations in a complementarity form at the element boundaries, which allows discontinuous approximations of the velocity.
Both elastic and inelastic impacts are treated.
However, it is not clear how state jumps in the tangential directions due to friction are treated in this work.
In addition, in the 3D case, they use a polyhedral approximation for the friction cone~\cite{Stewart1996a}.
In \cite{Shield2022}, only Gauss-Legendre implicit RK methods are regarded, which are known to have lower accuracy than L-stable methods such as Radau II-A for DAEs of index 3 \cite{Hairer1989}, e.g., during persistent contact motion phases.

\paragraph*{Contributions}
In this paper, we introduce the FESD-J method, which extends the FESD method \cite{Nurkanovic2022,Nurkanovic2022b,Nurkanovic2023} to rigid bodies with Coulomb friction and both elastic and inelastic impacts.
In particular, all the contacts can have different coefficients of restitution and friction.
Moreover, in the 3D case, we consider the exact nonlinear friction cone and do not make modeling approximations.
Our method can use any Runge-Kutta scheme and is thus more general than~\cite{Patel2019,Shield2022}.
Building on the ideas presented in~\cite{Shield2022}, we derive a complementarity form for the impact equations that can handle both normal and tangential state jumps, enabling a discontinuous approximation of the velocity state.
The impact equations are both on position and velocity level, and additional inequality constraints exclude spurious solutions in case of contact making.
We also establish cross-complementarity conditions that ensure accurate switch detection in the case of contact making and breaking, slip-stick transitions, or zero crossing of the tangential velocity in slip motions.
Finally, we derive step equilibration conditions that eliminate spurious degrees of freedom if no switches occur.
The presented method and its variations are implemented in the open-source toolbox NOSNOC~\cite{Nurkanovic2022}.
Numerical simulations and examples of optimal control problems are used to illustrate the effectiveness of the proposed method.
\paragraph*{Outline}
The paper is structured as follows.
Section \ref{sec:problem_formulation} introduces the dynamic model equations for rigid bodies subject to friction and impact.
Moreover, we propose a complementarity formulation for the impact conditions, which is convenient for event-based time discretizations.
In Section \ref{sec:fesd}, we introduce the novel \FESDj~method.
Section \ref{sec:fesd_ocp} discusses how to use the \FESDj\ method in the direct transcription of optimal control problems.
Section \ref{sec:examples} provides numerical examples.
The paper finishes with Section \ref{sec:conclusion}, which concludes the paper and discusses open problems.

\paragraph*{Notation}
All vector inequalities are to be understood element-wise, $\mathrm{diag}(x)\in \R^{n\times n}$ returns a diagonal matrix with $x \in \R^n$ containing the diagonal entries.
The concatenation of two column vectors $a\in \R^{n_a}$, $b\in \R^{n_b}$ is denoted by $(a,b)\coloneqq[a^\top,b^\top]^\top$, the concatenation of several column vectors is defined analogously.
The complementary conditions for two vectors  $a,b \in \R^{n}$ \sloppy{read as ${0\leq a \perp b\geq 0}$, where $a \perp b$ means $a^{\top}b =0$}.
The superscript $i$ denotes the contact index, e.g., $\Jni(q)$ is the contact normal of the $i$-th contact.
The subscripts "$\mathrm{n}$" and "$\mathrm{t}$" highlight if a variable or function is related to normal or tangential quantities. 
For the left and the right limits, we use the notation ${x(\ts^+)  = \lim\limits_{t\to \ts,\ t>\ts} x(t)}$ and ${x(\ts^-)  = \lim\limits_{t\to \ts,\  t<\ts}x(t)}$, respectively.
The identity matrix is denoted by $I \in \R^{n \times n}$ and a column vector with all ones is denoted by $e=(1,1,\dots,1) \in \R^n$, their dimension is clear from the context.
With $\mathbf{0}_{m,n} \in \R^{m\times n}$ we denote a matrix whose entries are all zeros.

\section{Problem formulation}\label{sec:problem_formulation}
In this section, we introduce the equations of motion for rigid bodies with impacts and Coulomb friction.
Subsection \ref{sec:cls} introduces the complementarity Lagrangian system.
In Subsection \ref{sec:eom}, we modify the basic equations of motion and in Subsection \ref{sec:friction_model} we discuss the complementarity formulation of the Coulomb friction model.
We finish this section by deriving the impulse equations in a convenient complementarity form.
Our modifications pave the way for deriving the \FESDj\ method.
\subsection{The complementarity Lagrangian system}\label{sec:cls}
There are many ways to represent nonsmooth rigid body systems, e.g., as second-order sweeping process \cite{Brogliato2016,Moreau1988,Moreau1999}, measure differential inclusion~\cite{Brogliato2016,Stewart2000,Stewart2011}, hybrid dynamical system~\cite{Johnson2016,Lunze2009}, differential variational inequalities of index 2~\cite{Stewart2000,Stewart2011} or Complementarity Lagrangian System (CLS)~\cite{Acary2008,Brogliato2016,Stewart2011}.
In this paper, we focus on the CLS representation, which is the usual choice in the development of time-discretization methods.
The CLS model equations with friction read as:
\begin{subequations}\label{eq:cls}
	\begin{align}
		\dot{q} &= v, \label{eq:cls_pos}\\
		M(q)\dot{v} &= f_\mathrm{v}(q,v) + B_u(q)u  +  \sum_{i=1}^{\ncontacts} (\Jni(q) \zni+ \Jti(q)\zti), \label{eq:cls_forces}\\
		0 &\leq \zni \perp \fci(q) \geq 0, &i=1,\ldots,\ncontacts, \label{eq:cls_complementarity}\\
			0&=\Jni (q(\ts))^\top (v(\ts^+)+\CoRi v(\ts^-)), \;\
			\textrm{if} \ \fci(q(\ts)) = 0 \; \textrm{and} \; \Jni(q(\ts))^\top v(\ts^-)<0, &i=1,\ldots,\ncontacts,\label{eq:cls_state_jump}\\
		&\zti  \in\arg\min_{{\zttildei} \in \R^{\nt}} \quad  -v^\top \Jti(q){\zttildei}
		\
		\quad	\textrm{s.t.} \quad  \| {\zttildei} \|_2 \leq \mu^{i} \zni, &i=1,\ldots,\ncontacts. \label{eq:cls_friction_disipiation}
	\end{align}
\end{subequations}
For ease of notation, we omitted the time dependencies in the differential and algebraic states.
The matrix $M(q)\in \R^{n_q\times n_q}$ is the inertia matrix and is assumed to be symmetric positive definite. 
The differential state are the generalized coordinates $q(t)\in\R^{n_q}$ and generalized velocities $v(t)\in\R^{n_q}$.
The function $f_v(q,v)\in \R^{n_q}$ collects gravitational, Coriolis, and all other external forces. 
The matrix $B(q) \in \R^{n_q \times n_u}$ is the control input mapping and we assume that the generalized control force $u(t) \in \R^{n_u}$ is known, e.g., obtained by solving an Optimal Control Problem (OCP).

We say that we have $\ncontacts$ possible frictional contacts indexed by $i$, which are modeled by the unilateral constraints $\fci(q) \geq 0$, $i =1,\ldots,\ncontacts$.
The smooth scalar functions $\fci(q)$ represent the signed distances.
Consequently, we have the generalized normal contact forces, $\zni \in \R$, and tangential contact forces, $\zti \in \R^{\nt}$, acting on the rigid body.
For planar contacts we have $\nt= 1$ and for 3D contacts $\nt =2$.
All algebraic variables are grouped into the vectors $\lambda = (\zn,\zt) \in \R^{\ncontacts(1+\nt)}$, where
$\zn =(\lambda_{\mathrm{n}}^{1},\ldots,\lambda_{\mathrm{n}}^{\ncontacts})$ and
$\zt =(\lambda_{\mathrm{t}}^{1},\ldots,\lambda_{\mathrm{t}}^{\ncontacts})$.

The complementarity constraint \eqref{eq:cls_complementarity} expresses for every contact that: either the body is not in contact ($\fci(q)>0$) and there is no normal contact force $(\zni =0)$, or the body is in contact ($f_{c,i}(q) =0$) and a nonnegative normal contact force ($\zni\geq0$, i.e., there is no adhesion) acts along the surface contact normal $\Jni(q) \coloneqq \nabla_q \fci(q) \in \R^{n_q}$.
Whenever a contact becomes active with a respective negative normal velocity $\Jni(q)^\top v <0$, then a state jump must occur to avoid interpenetration.
This is expressed with \eqref{eq:cls_state_jump}, where the first part is Newton's restitution law, which expresses the post-impact normal velocity $\Jni(q(\ts))^\top v(\tsp) = \CoRi \Jni(q(\ts))^\top v(\tsn) \geq0$, as a function of the pre-impact velocity.
The scalar $\CoRi\in [0,1]$ is the \textit{coefficient of restitution} of the $i$-th contact.
For $\CoRi = 0$ we speak of inelastic impacts and for $\CoRi>0$ of (partially) elastic impacts.

The last part of the CLS is the Coulomb friction model expressed via the convex optimization problem~\eqref{eq:cls_friction_disipiation}. 
The friction force of the $i$-th contact acts in the tangent space at the contact point, which is spanned by the tangent Jacobian $\Jti(q) \in \R^{n_q \times \nt}$. 
It models the \textit{maximum dissipation principle}, which expresses that the friction force $\zti$ is chosen such that it maximizes the energy dissipation.
This friction law has the following features:
\begin{itemize}
	\item the friction force has the opposite direction to the tangential slip direction $\Jti(q)^\top v$
	\item the maximal magnitude of the friction force $\zti$ is the product of $\mu^{i}$, the \textit{coefficient of friction} of the $i$-th contact, and the normal contact force $\zni$.
\end{itemize}
Note that if we have nonzero tangential velocity, i.e., $\Jti(q)^\top v \neq 0$, then $\|\zti\|_2  = \mu \zni$ (slip motion phase) and for $\Jti(q)^\top v = 0$, we have $\|\zti\|_2  \leq \mu \zni$ (stick phase).
The friction model is another source of nonsmoothness.
Note that extensions to anisotropic friction are straightforward if we replace the Euclidian norm in~\eqref{eq:cls_friction_disipiation} with an appropriately weighted norm.

The set of all possible contact forces at the $i$-th contact is called the friction cone and is defined as
\begin{align*}
	\mathrm{FC}_{i}(q) = \{ \Jni(q)\zni + \Jti(q) \zti \mid \zni \geq 0, \|\zti \|_2 \leq \mu \zni  \}.
\end{align*}
The total friction cone is the sum of all friction cones generated by each contact:
\begin{align*}
	\mathrm{FC}(q) = \sum_{i \in \{j \mid f_{c,j}(q) =0\}} \mathrm{FC}_{i}(q).
\end{align*}

We assume that the functions $M(q), f_v(q,v), B(q), f_c(q), \Jni(q) $ and $\Jti(q)$ are at least twice continuously differentiable.
We proceed by reformulating some parts of \eqref{eq:cls} to obtain a formulation more suitable for discretization via the \FESDj\ method presented below.
\subsection{Modified equations of motion}\label{sec:eom}
For ease of notation, we rewrite the velocity dynamics equation as
\begin{align}\label{eq:velocity_dynamics}
\dot{v} &= F_v(q,v,u,\lambda) \coloneqq M(q)^{-1}(f_\mathrm{v}(q,v) + B_u(q)u  +  \sum_{i=1}^{\ncontacts} (\Jni(q) \zni+ \Jti(q)\zti))
\end{align}
In general, it might be costly to symbolically invert the matrix $M(q)$.
Alternatively, one can introduce a \textit{lifting variable} $z_v\in \R^{n_q}$ and impose the equations: 
\begin{align*}
\dot {v} &= z_v,\\
M(q)z_v &= f_\mathrm{v}(q,v) + B_u(q)~u  +  \sum_{i=1}^{\ncontacts} (\Jni(q) \zni+ \Jti(q)\zti),
\end{align*}
where the inverse is computed during the numerical solution of this algebraic equation within the time-discretization method.

\subsection{Friction model}\label{sec:friction_model}
We discuss the complementarity formulations of the nonlinear friction cone model in \eqref{eq:cls_friction_disipiation} and a commonly used polyhedral approximation~\cite{Stewart2000,Stewart1996a}.
In the planar cases, they are equivalent and in the 3D case, the polyhedral cone is an approximation that can be made arbitrarily accurate.
\subsubsection{Nonlinear friction cone}
Let us introduce the shorthand for the tangential velocity at the $i$-th contact point $\vti \coloneqq {\Jti}^\top v \in \R^{\nt}$.
For a given normal contact force $\zni$, the solution map of \eqref{eq:cls_friction_disipiation} reads as:
\begin{align}\label{eq:friction_slution_map}
	\zti &\in  \begin{cases}
		\{- \mu^{i} \zni  \frac{ \vti }{\| \vti  \|_2} \}, & \textrm{if}\; \| \vti \|_2 > 0,\\
		\{ {\zttilde} \mid \|\zttilde\|_2 \leq \mu^{i} \zni \}, & \textrm{if}\; \| \vti \|_2 = 0.
	\end{cases} 
\end{align}
Note that in the planar case, this reduces to the well-known expression $\zti \in -\mu^{i} \zni \sign(\vti)$. 
For differentiability, we replace the inequality constraint in \eqref{eq:cls_friction_disipiation} by the equivalent constraint
\begin{align}\label{eq:squared_conict_cstr}
& \| {\zttildei} \|_2^2 \leq (\mu^{i} \zni)^2.
\end{align}
The KKT conditions of the modified system (\eqref{eq:cls_friction_disipiation}, now with the inequality constraint \eqref{eq:squared_conict_cstr}) for $i = 1,\ldots,\ncontacts$ read as:
\begin{subequations}\label{eq:friction_equations_3D}
\begin{align}
	&0= -\Jti(q)^\top v -2\gamma^{i} \zti,\\
	&0 = \beta^i - (\mu^{i} \zni)^2 + \| {\zti} \|_2^2, \label{eq:friction_equations_3D_lift}\\
	&0 \leq \gamma^{i} \perp \beta^i\geq 0 \label{eq:friction_equations_3D_comp}.
\end{align}
\end{subequations}
where $\gamma^{i} \in \R$ is the Lagrange multiplier for \eqref{eq:squared_conict_cstr} and $\beta^i$ is an auxiliary lifting variable used to remove the nonlinearity from the complementarity condition~\eqref{eq:friction_equations_3D_comp}.

If the $i$-th contact is active, with $\zni >0$, and no further impacts occur, then  $\vti = \Jti(q)^\top v$ is continuous in time.
Suppose that is $\zni$ continuous during persistent contact.
By manipulating \eqref{eq:friction_equations_3D}, we can obtain that $\gamma^{i} = \frac{\|\vti\|}{2\mu^{i} \zni}$, which implies 
that $\gamma^{i}$ is also a continuous function of time.
To achieve switch detection in all cases, we will need to look at the positive and negative parts of the tangential velocity components, cf. Section \ref{sec:fesd_cross_comp}.
We can obtain them by introducing the two slack variables $\xi^{i}_{+},\xi^{i}_{-} \in \R^{\nt}$ and augmenting \eqref{eq:friction_equations_3D} by: 
\begin{subequations}\label{eq:friction_equations_3D_augmented}
\begin{align}
	&0 = \Jti(q)^\top v -\xi^{i}_{+}+\xi^{i}_{-},\\
	&0 \leq \xi^{i}_{+} \perp \xi^{i}_{-}  \geq 0.
\end{align}
\end{subequations}
It can be seen that $\xi^{i}_{+}  = \max(0, \Jti(q)^\top v)$ and $\xi^{i}_{-}  = \max(0, -\Jti(q)^\top v)$.
If no impact occurs, the tangential velocity $\Jti(q)^\top v$ is a continuous but not necessarily smooth function of time.
Therefore, the variables $\xi^{i}_{+}$ and $\xi^{i}_{-}$ are also continuous functions of time.
These continuity properties in the formulations above are essential for the derivation of the \FESDj\ discretization for CLS and its exact switch detection capabilities.
\subsubsection{Polyhedral friction cone}
The nonlinear 3D friction cone is often approximated by a polyhedral cone~\cite{Stewart1996a}, which yields a linear programming friction model.
However, the polyhedral friction cone can introduce drifting if the tangential velocity is not aligned with some of the vectors spanning the tangent plane~\cite{Stewart1996}.
Regard the matrix $\Dti(q) \in \R^{n_q \times n_d}$, such that for every column $d^{i}_k(q)$ of $\Dti(q)$ there exist another column vector $d^{i}_{j}(q)$ such that $d^{i}_{k}(q)=-d^{i}_{j}(q)$. 
An example choice is $\Dti(q) =\begin{bmatrix}	\Jti(q)& -\Jti(q) \end{bmatrix}$.
For a fixed $\zni>0$, the convex combination of the columns of $\Dti(q)$ defines a polyhedron that approximates the circle \eqref{eq:squared_conict_cstr}, for an illustration cf. \cite[Figure 2.2]{Stewart2000}.
Therefore, the friction force is approximated by
\begin{align*}
\sum_{i=1}^{\ncontacts} \Dti(q) \zdi \approx \sum_{i=1}^{\ncontacts}\Jti(q)\zt.
\end{align*}
The convex problem \eqref{eq:cls_friction_disipiation} in the CLS is replaced by the linear program:
\begin{subequations}\label{eq:friction_cone_polyhedral_lp}
	\begin{align}
		\min_{\zdi \in \R^{n_d}} \ &-v^\top \Dti(q) \zdi \\
		\quad \mathrm{s.t.} \quad & e^\top \zdi \leq \mu^i \zni,\\
		&\zdi \geq0.
	\end{align}
\end{subequations}
Using its KKT conditions and expressing the Lagrange multipliers for the second inequality from the stationarity conditions, we obtain for all $i=1,\ldots,\ncontacts$ the following set of complementarity conditions:
\begin{subequations}\label{eq:friction_cone_polyhedral_comp}
	\begin{align}
		&0\leq \zdi \perp {\Dti(q)}^\top v+ \gamma_d^{i} e \geq 0,\\
		&0\leq \gamma_d^{i}\perp \mu^i \zni - e^\top \zdi\geq 0,
	\end{align}
\end{subequations}
where $\gamma_d \in \R$ approximates the magnitude of the tangential contact velocity and is, if no impact occurs a continuous function of time~\cite{Stewart2000}, whereas the functions $\zdi$ do not have in general continuity properties.
The polyhedral approximation \eqref{eq:friction_cone_polyhedral_lp} is in the planar case exact.
Note that in the 3D case, for a single contact, if we have at least four columns in $\Dti(q)$, Eq.~\eqref{eq:friction_cone_polyhedral_comp} has more complementarity constraints than the augmented nonlinear model \eqref{eq:friction_equations_3D}-\eqref{eq:friction_equations_3D_augmented}.
In the remainder of this paper, we focus on the nonlinear friction cone. 
The developments for the polyhedral case follow similar lines. 
Implementations of both variants are available in NOSNOC. 
\subsection{Impulse equations and state jump laws}\label{sec:impulse_equations}
The point-wise algebraic law $\eqref{eq:cls_state_jump}$ complicates the numerical treatment of the CLS.
Let us regard a single contact and suppose for a moment that there is no friction, i.e., $\zt =0$ and we omit Eq. \eqref{eq:cls_friction_disipiation} in the CLS model. 
At the time of impact, the normal contact force $\zni$ is a Dirac "function".
If we would know that at time $\ts$ an impact would occur at the $i$-th contact, to obtain the post-impact velocity $v(\ts^+)$, we can solve the point-wise \textit{impulse} equations \cite{Brogliato2016}:
\begin{subequations}\label{eq:impulse_eq_plain}
\begin{align}
	&M(q(\ts))(v(\tsp)-v(\tsn)) = \Jni(q(\ts))\Zni,\\
	&\Jni(q(\ts))^\top (v(\tsp) + \CoRi v(\tsn)) =0,
\end{align}
\end{subequations}
where $\Zni$ is the normal contact impulse, i.e., the integral of the normal contact force at the impact 
$\Zni = \lim_{\epsilon\to0} \int_{\ts-\epsilon}^{\ts+\epsilon}\zni(t) \dd t$.
In the sequel, the function $q(\cdot)$ is always evaluated at $\ts$ and we omit this dependency for brevity.

However, in practice, we usually do not know when and how many impacts will occur.
We propose to rewrite the impulse equations in complementarity form both at position and velocity level:
\begin{subequations}\label{eq:impulse_eq_comp}
	\begin{align}
		&M(q)(v(\tsp)-v(\tsn)) = \Jni(q)\Zni,\\
		&0\leq \Zni \perp \fci(q)\geq 0 \label{eq:impulse_eq_comp_pos},\\
		&0= \Zni  \Jni(q)^\top (v(\tsp) + \CoRi v(\tsn)) \label{eq:impulse_eq_comp_vel}.
	\end{align}
\end{subequations}
Note that we do not impose nonnegativity constraints on $\Jni(q)^\top(v(\tsp) + \CoRi v(\tsn))$, which is usually done when writing the impulse equation on velocity level~\cite{Brogliato2016}.
We do this to be able to use \eqref{eq:impulse_eq_comp} even if no contacts happen, i.e., $\fci(q)>0$. 
In this case, it is also possible that $\Jni(q)^\top (v(\tsp) + \CoRi v(\tsn))<0$.
This system of equations \eqref{eq:impulse_eq_comp} is over-determined, but still yields well-defined solutions and efficiently encodes different impacts, as we discuss next.
On the one hand, if $\fci(q) >0$ (no contact), we have that $\Zni =0$, and since the matrix $M(q)$ is full rank, we have that the velocity stays continuous, i.e., $v(\tsp)= v(\tsn)$. 
Since, $\Zni = 0$, equation \eqref{eq:impulse_eq_comp_vel} is implicitly satisfied it does not affect the solution.
On the other hand, for $\fci(q) = 0$, we have from \eqref{eq:impulse_eq_comp_pos} that $\Zni \geq 0$.
In this case, equation \eqref{eq:impulse_eq_comp_vel} helps us to determine the value of $\Zni$, and we obtain the following reduced linear system:
\begin{align*}
		&M(q)(v(\tsp)-v(\tsn)) = \Jni(q)\Zni,\\
		&0= \Jni(q)^\top (v(\tsp) + \CoRi v(\tsn)),
\end{align*}
that permits two solutions. 
First if, $\Zni >0$, then we can compute that
\begin{align*}
\Zni = - \big(\Jni(q)^\top M(q)^{-1} \Jni(q)\big)^{-1} (1+\CoRi)\Jni(q)^\top v(\tsn),\\
\Jni(q)^\top v(\tsp) =  \Jni(q)^\top M(q)^{-1} \Jni(q)\Zni +\Jni(q)^\top v(\tsn) \geq 0.
\end{align*}
Second, if $\Zni = 0$, we obtain that $v(\tsp)= v(\tsn)$. 
In an impacting scenario thus would imply that $\Jni(q)^\top v(\tsp) = \Jni(q)^\top v(\tsn) <0$. 
However, this would lead to the violation of $\fci(q(t)) \geq 0 $ for $t> \tsp$. 
Consequently, only the first solution is feasible for the overall CLS.

A similar formulation that aggregates the position and velocity constraints with the help of slack variables is provided in \cite{Shield2022}.
For the reader's convenience, we restate them in ~\ref{app:shield_equations}.

We proceed with extending these equations to the frictional case.
The normal contact impulse $\Zni$  results in a tangential impulse $\Zti \in \R^{\nt}$.
Following the contact model used in \cite{Anitescu1997}, we require that the maximum dissipation principle be satisfied for the contact impulses at the post-impact tangential velocity $\vti(\tsp)$. 
Using \eqref{eq:friction_equations_3D}, and combining it with \eqref{eq:impulse_eq_comp}, we obtain the full impact model for the $i$-th contact as:
\begin{align*}
		&M(q)(v(\tsp)-v(\tsn)) = \Jni(q)\Zni+\Jti(q)\Zti,\\
		&0\leq \Zni \perp \fci(q)  \geq0,\\
		&0= \Zni  \Jni(q)^\top (v(\tsp) + \CoRi v(\tsn)),\\
		&0=-\vti(\tsp)  - 2\Gamma^{i} \Zti,\\
		&0 = B^{i} - (\mu^{i} \Zni)^2 + \|\Zti\|_2^2,\\
		&0 \leq \Gamma^{i} \perp B^i \geq 0,
\end{align*}
where $\Gamma^{i}$ is the Lagrange multiplier corresponding to $\|\Zti\|_2 \leq (\mu^{i} \Zni)^2 $ and $B^i$ is a lifting variable.
In our notation, all \textit{impulse-related} quantities are denoted by uppercase letters. 
In the case of multiple simulations impacts, the equations are simply aggregated and we obtain the impact model:
\begin{subequations}\label{eq:impulse_eq_comp_friction_all}
	\begin{align}
		&M(q)(v(\tsp)-v(\tsn)) = \sum_{i=1}^{\ncontacts} \Jni(q)\Zni+\Jti(q)\Zti,\\
		&0\leq \Zni \perp \fci(q) \geq0,\ \quad &i = 1,\ldots,\ncontacts, \\
		&0= \Zni \Jni(q)^\top (v(\tsp) + \CoRi v(\tsn)),\quad &i = 1,\ldots,\ncontacts, \\
		&0=-\vti(\tsp)  - 2\Gamma^{i} \Zti,\ \quad &i = 1,\ldots,\ncontacts, \\
		&0= B^i - (\mu^{i} \Zni)^2+\|\Zti\|_2^2, \quad &i = 1,\ldots,\ncontacts, \\
		&0 \leq \Gamma^{i} \perp B^i  \geq 0, &\quad i = 1,\ldots,\ncontacts. 
	\end{align}
\end{subequations}
Note that there might be different simultaneous contact models \cite{Nguyen2018} and which is used is a modeling decision.
We take this model, which is commonly used in time-stepping methods and has proven to provide physically realistic results~\cite{Anitescu2006,Anitescu2004,Anitescu1997,Nguyen2018,Stewart2000}.
In a concrete application, one should always verify the physical soundness of the produced results.

To summarize, during impacts, the position stays continuous $q(\tsp) = q(\tsn)$, and the velocities are updated in the according directions, since $\Zn \neq 0$ and $\Zti \neq0$ (if $\mu^{i} >0$ and $\|\vti\|_2 \neq 0$). 
On the other hand, if $\fc(q)>0$, then $\Zn=0$, $\Zt= 0$, and the velocity stays continuous as well $v(\tsp) = v(\tsn)$.
The formulation \eqref{eq:impulse_eq_comp_friction_all} helps to encode the continuity and state jump equations in one shot, for both elastic and inelastic impacts.

\section{Finite Elements with Switch and Jump Detection (\FESDj)}\label{sec:fesd}
In this section, we derive the FESD-J method for \eqref{eq:cls}.
We discretize the equations of motions and the complementarity conditions with a standard Runge-Kutta method.
Additionally, we let the integration step size $h_n$ be degrees of freedom.
Similarly to \cite{Baumrucker2009,Nurkanovic2022}, we introduce \textit{cross complementarity} conditions that make active set changes, switches, and impacts possible only at the boundaries of the finite element (integration step).
This will lead to exact switch detection. 
However, if no switches occur, for a fixed active set Eq. \eqref{eq:cls} reduces to a smooth ODE or differential algebraic equations. 
Consequently, the step sizes $h_n$ are spurious degrees of freedom in this case.
By adapting the ideas from \cite{Nurkanovic2022}, we introduce \textit{step equilibration} conditions, which resolve this degeneracy.
As in any event-based method, we assume that there are sufficiently many finite elements $\NFE$ to be able to detect all switches.

\FESDj\ should detect several different kinds of discontinuities: 
1) discontinuities due to contact making and breaking, which usually introduce discontinuities in the velocities,
2) discontinuities due to friction during persistent contact: transitions from slip to stick, from stick to slip, or zero crossing of tangential velocities in slip motion, which all may lead to discontinuities in the friction forces and kinks in the velocity.

\subsection{Runge-Kutta discretization}
We regard a single control interval $[0,T]$ with a constant known control input $\hat{u} \in \R^{n_u}$, i.e., we set $u(t) = \hat{u}$ for $ t\in [0,T]$.
Section \ref{sec:fesd_ocp} will treat the discretization of optimal control problems with multiple control intervals.
The time interval $[0,T]$ is divided into $\NFE$ \textit{finite elements} $[t_{n},t_{n+1}]$, with the grid points $0= t_0 < t_1 < \ldots <t_{\NFE} = T$.
On each of the finite elements we regard an $\Nstg$-stage Runge-Kutta method which is parameterized by its Butcher tableau entries $a_{i,j} ,b_i$ and $c_i$ with $i,j\in\{1,\ldots,\Nstg\}$~\cite{Hairer1991}.
The integration step sizes are defined as the finite element lengths, i.e., $h_{n} = t_{n+1} - t_{n},\; n = 0, \ldots,\NFE-1$.

\begin{figure}[t]
	\centering
	{\includegraphics[width=0.75\textwidth]{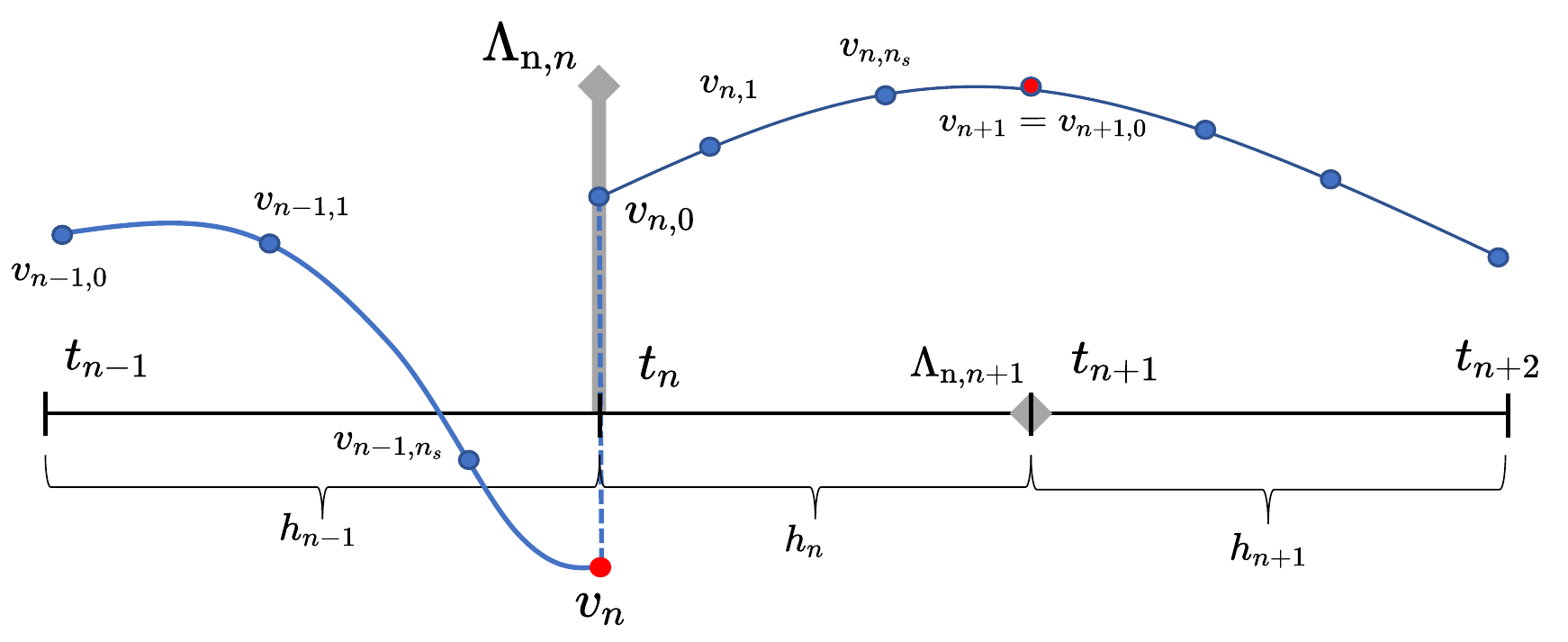}}
	\caption{Illustration of the \FESDj\ approximation for the velocity state. 
	At $t_n$ an impact occurs and the impulse approximation $\Znn$ leads to a discontinuity in the velocity, i.e., $v_{n} \neq v_{n,0}$.
	At time $t_{n+1}$, there is no impact, $\Lambda_{\mathrm{n},n+1} =0$, and the velocity is continuous, , i.e., $v_{n+1} = v_{n+1,0}$.}
	\label{fig:fesd_illustration}
\end{figure}

The approximations of the differential states at the grid points $t_n$ are denoted by $q_n \approx q(t_n)$, $v_n \approx v(t_n)$.
Their values at the RK stage points $t_{n,i} \coloneqq t_n + c_i h_n,\; i = 1,\ldots, \Nstg$, are denoted by $q_{n,i}$, $v_{n,i}$.
Likewise, the stage values of the algebraic variables are denoted by $\lambda_{\mathrm{n},n,j}^{i}$,$\lambda_{\mathrm{t},n,j}^{i}$, $\gamma_{n,j}^{i}$, $\beta_{n,j}^{i}$,
$\xi_{+,n,j}^{i}$ and $\xi_{-,n,j}^{i}$.
Moreover, we denote the differential state approximations at the left boundary point of a finite element by $q_{n,0}$ and $v_{n,0}$, which corresponds to $q(t_n^+)$ and $v(t_n^+)$. 
In the case of state jumps we have at $t_n$ that $v(t_n^-) \neq v(t_n^+)$. 
Similarly, in the time-discretization derived below, we may have $v_{n} \neq v_{n,0}$.
Figure~\ref{fig:fesd_illustration} provides a schematic illustration of our notation conventions.
For ease of notation, we assume that $c_{\Nstg} = 1$ and comment on how to extend to the case when $c_{\Nstg} \neq 1$, e.g. in Gauss-Legendre implicit RK schemes, later in Remark~\ref{remark:right_boundary_point}.

Given the initial value $s_0 \coloneqq (q_0,v_0)$, we first discretize the differential equations for the position \eqref{eq:cls_pos} and velocity \eqref{eq:cls_forces} for all $n = 0,\ldots, \NFE-1$:
\begin{subequations}\label{eq:irk_ode_int}
	\begin{align}
		&0= q_{n,j} - (q_{n,0} + h_n \sum_{k=1}^{\Nstg} a_{j,k} v_{n,k}),\quad &j =1,\ldots,\Nstg,\\
		&0= v_{n,j} - (v_{n,0} + h_n \sum_{k=1}^{\Nstg} a_{j,k} F_v(v_{n,k},q_{n,k},\lambda_{n,k},\hat{u})),\quad &j =1,\ldots,\Nstg,\\
		&0= q_{n+1} - (q_{n,0} + h_n \sum_{k =1}^{\Nstg} b_k  v_{n,k}), \label{eq:irk_q_out}\\
		&0= v_{n+1} - (q_{n,0} + h_n \sum_{k =1}^{\Nstg} b_k F(v_{n,k},q_{n,k},\lambda_{n,k}, \hat{u})), \label{eq:irk_v_out}
	\end{align}
\end{subequations}
where $\lambda_{n,k} = (\lambda_{\mathrm{n},n,k},\lambda_{\mathrm{t},n,k}) \in \R^{\ncontacts(1+n_\mathrm{t})}$ collects the stage values for the normal and tangential contact force approximations.
The algebraic contact conditions \eqref{eq:cls_complementarity} and friction model equations \eqref{eq:friction_equations_3D} are evaluated for $n = 0,\ldots,\NFE-1$, at the RK stage points:
\begin{subequations}\label{eq:irk_ode_alg}
	\begin{align}
		&0\leq \fci(q_{n,j}) \perp \lambda_{\mathrm{n},n,j}^{i} \geq 0,\quad &i =1,\ldots,\ncontacts,\ j =1,\ldots,\Nstg, \label{eq:irk_ode_alg_contact_comp}\\
		&0= -\Jti(q_{n,j})^\top v_{n,j} +2\gamma^{i}_{n,j} \lambda^{i}_{\mathrm{t},n,j},\quad &i =1,\ldots,\ncontacts,\ j =1,\ldots,\Nstg,\\
		&0 =\beta^{i}_{n,j} -  (\mu^{i} \lambda_{\mathrm{n},n,j})^2 + \| \lambda_{\mathrm{t},n,j} \|_2^2,\quad &i =1,\ldots,\ncontacts,\ j =1,\ldots,\Nstg,\\
		&0 \leq \gamma^{i}_{n,j} \perp \beta^{i}_{n,j} \geq 0,\quad &i =1,\ldots,\ncontacts,\ j =1,\ldots,\Nstg. \label{eq:irk_ode_alg_friction1}
	\end{align}
\end{subequations}
Additionally, we discretized the auxiliary conditions \eqref{eq:friction_equations_3D_augmented}
for all $n=0,\ldots,\NFE-1$, at all RK stage points $j = 1,\ldots, \Nstg$, and at the left boundary point of every finite element (we assign the index $j=0$ for these variables):
\begin{subequations}\label{eq:friction_pos_neg_discrete}
	\begin{align}
	&0 =\Jti(q_{n,j})^\top v_{n,j} -\xi^{i}_{+,n,j} + \xi^{i}_{-,n,j},\quad &i =1,\ldots,\ncontacts,\ j =0,\ldots,\Nstg \label{eq:friction_pos_neg_def},\\
	&0 \leq \xi^{i}_{+,n,j} \perp \xi^{i}_{-,n,j} \geq 0,\quad &i =1,\ldots,\ncontacts,\ j =0,\ldots,\Nstg. \label{eq:friction_pos_neg_comp}
	\end{align}
\end{subequations}
Including the left boundary point is crucial for switch detection, as we show below.

\begin{remark}(Extension for $c_{\Nstg} \neq 1$)
\label{remark:right_boundary_point}
In the cross complementarity conditions introduced in Section \ref{sec:fesd_cross_comp}, we exploit the continuity properties of some of the algebraic variables.
For the differential and some algebraic states we need the left and right boundary points of a finite element.
If $c_{\Nstg} \neq 1$, then the right boundary point is not an RK stage point.
We only must extend~\eqref{eq:friction_pos_neg_discrete} by adding the constraints for all $n=0,\ldots,\NFE-1$:
\begin{align*}
	&0 =\Jti(q_{n+1})^\top v_{n+1} -\xi^{i}_{+,n+1} + \xi^{i}_{-,n+1}, &i =1,\ldots,\ncontacts,\\
	&0 \leq \xi^{i}_{+,n+1} \perp \xi^{i}_{-,n+1} \geq 0, &i =1,\ldots,\ncontacts.
\end{align*}
The variables $v_{n+1}$ and $ q_{n+1} $ are naturally obtained from the IRK discretization via \eqref{eq:irk_q_out} and \eqref{eq:irk_v_out}, while the algebraic variables $\xi^{i}_{+,n+1}$ and $ \xi^{i}_{-,n+1}$ need to be additionally defined for the right boundary point of the $n$-th finite element.
\end{remark}

As the step sizes are degrees of freedom, we require that their sum equals the total interval length:
\begin{align}\label{eq:sum_h}
	&0= T - \sum_{n=0}^{\NFE-1} h_n.
\end{align}

\subsection{Continuity and state jumps}
We proceed with linking the right boundary points of a finite element $q_{n}$, $v_{n}$ with the left boundary points $q_{n,0}$,$v_{n,0}$, of the following finite element.
The generalized positions are continuous functions of time and we impose for their discrete-time counterparts for all $n = 0,\ldots,\NFE$:
\begin{align}\label{eq:continuity_q}
	&0=q_{n,0} - q_{n}.
\end{align}
The junction conditions for the velocities are more elaborate as they can be continuous or have jumps in both normal and tangential directions.
We express all the possible combinations at once with the discrete-time version of \eqref{eq:impulse_eq_comp_friction_all}.
Recall that $v_n$ corresponds to the pre-impact velocity $v(t_n^-)$ and $v_{n,0}$ to the post-impact velocity $v(t_n^+)$.
Thus, for all $n = 0,\ldots, \NFE-1$ we have:
\begin{subequations}\label{eq:impulse_eq_comp_friction_all_discretized}
	\begin{align}
		&0=M(q_n)(v_{n,0}-v_{n}) -\sum_{i=1}^{\ncontacts} \Big(\Jni(q_{n})\Znin+\Jti(q_{n})\Ztin\Big), \\
		&0\leq \Znn^{i} \perp \fci(q_{n})\geq0,								\quad & i = 1,\ldots,\ncontacts,\label{eq:position_switch_detection} \\
		&0=\Znn^{i} \Jni(q_{n})^\top (v_{n,0} + \CoRi v_{n}), 				\quad & i = 1,\ldots,\ncontacts, \\
		&0=-\Jti(q_{n})^\top v_{n,0}  - 2\Gamma_{n}^{i} \Ztin,\ 			\quad & i = 1,\ldots,\ncontacts, \label{eq:velocity_switch_detection}\\
		&0 = B_n^{i} - (\mu^{i} \Znin)^2 + \|\Ztin\|_2,\ 			\quad & i = 1,\ldots,\ncontacts,\\ 
		&0 \leq \Gamma^{i}_n \perp B_n^{i} \geq 0,\ 						\quad & i = 1,\ldots,\ncontacts.
	\end{align}
\end{subequations}
Observer that if for the initial conditions for $n=0$, it holds that some $\fci(q_0) =0$ and $\Jni(q_0)^\top v_0<0$ a state jump will occur and reinitialize the velocity $v_{0,0}$ to a feasible value.
In the sequel, we introduce conditions that will not allow active set changes in the complementarity conditions (i.e., switches) on the RK-stage points within the finite element, cf. Figure~\ref{fig:fesd_illustration}.
Consequently, the jumps can happen only at the grid points $t_n$.
Therefore, impulse equations are evaluated only at the grid points $t_n, n = 0,\ldots,\NFE$.
Note that if $\fc(q_n) > 0$, then $\Znn = 0$, $\Ztn = 0$, and \eqref{eq:impulse_eq_comp_friction_all_discretized} implies the continuity condition $v_{n,0}=v_{n}$.

\begin{remark}{(Excluding continuity conditions in case of impacts)}
Recall from the discussion at the beginning of Section~\ref{sec:impulse_equations}, that the constraints $f_c(q(t)) \geq 0$ prevent post-impact velocities $v(t_n^+) = v(t_n^-)$, even if we do not explicitly impose the constraints $\Jni(q(t_n))^\top v(t_n^+) + \CoRi \Jni(q(t_n))^\top v(t_n^-) \geq 0$.
In the discrete-time setting, the non-negativity of the gap function is only imposed at the RK stage points, and at the finite elements boundaries, cf. Eq. \eqref{eq:irk_ode_alg_contact_comp} and Eq. \eqref{eq:position_switch_detection}, respectively.
However, this might not always be sufficient to exclude solutions with $\Zni = 0$ and very large $\zni >0$, which may lead to the satisfaction of the non-negativity of gap function constraints.
Therefore, we propose to add the following inequality constraints:
\begin{align}
	&\fci(q_{n,0} + \epsilon_{\mathrm{c}} h_n v_{n,0} ) \geq0, &n = 0,\ldots,\NFE-1,\ i = 1,\ldots,\ncontacts,
\end{align}
where $\epsilon_\mathrm{c}$ is a small positive number, e.g. $\epsilon_\mathrm{c} = 10^{-3}$.
The point $q_{n,0} + \epsilon_{\mathrm{c}} h_n v_{n,0}$ corresponds to an explicit Euler step of length $\epsilon_\mathrm{c} h_n $ and the non-negativity of the gap functions is imposed there, i.e., shortly after the impact.
In case of impacts, this excludes solutions with $\Zni =0$ and $\zni \gg 0$.
It can be seen that this constraint is always satisfied for the continuous time variant of the system and only excludes spurious solutions.
\end{remark}

We introduce a more compact notation for the conditions and variables introduced so far.
All position and velocity variables (stage and boundary points) are collected into the vectors
$\mathbf{q} = (q_0,q_{0,0},\ldots,q_{0,\Nstg},q_{1},\ldots,q_{\NFE-1,0},\ldots,q_{\NFE})$,
$\mathbf{v} = (v_0,v_{0,0},\ldots,v_{0,\Nstg},v_{1},\ldots,v_{\NFE-1,0},\ldots,v_{\NFE})$, respectively, and we define $\mathbf{x} = (\mathbf{q},\mathbf{v})$.
We collect all step sizes into $\mathbf{h} = (h_0,\ldots,h_{\NFE-1})$.
The stage values of all algebraic variables $\lambda_{\mathrm{n},n,j}^{i}$,$\lambda_{\mathrm{t},n,j}^{i}$, $\gamma_{n,j}^{i}$, $\beta_{n,j}^{i}$,
$\xi_{+,n,j}^{i}$, $\xi_{-,n,j}^{i}$, for all $ i=1,\ldots,\ncontacts, n = 0,\ldots,\NFE-1$ and $j = 1,\ldots,\Nstg$,
and
all algebraic variables appearing in the discrete-time impulse equations \eqref{eq:impulse_eq_comp_friction_all_discretized}, namely
$\Lambda_{\mathrm{n},n}^{i}$, $\Lambda_{\mathrm{t},n}^{i}$, $\Gamma_{n}^{i}$, $B_n^i$, for all $ i=1,\ldots,\ncontacts$ and $n = 0,\ldots,\NFE$ and $j = 1,\ldots,\Nstg$ are collected in the vector $\mathbf{z}$.
By noting that complementarity conditions $ 0 \leq a \perp b \geq 0$ can be equivalently written as a nonsmooth equations via C-functions $\psi(a,b) =0$, we collect the equations \eqref{eq:irk_ode_int}, \eqref{eq:irk_ode_alg}, \eqref{eq:friction_pos_neg_discrete}, \eqref{eq:continuity_q} and \eqref{eq:impulse_eq_comp_friction_all_discretized} into the equation:
\begin{align}\label{eq:rk_equations}
	&0 = G_{\mathrm{rk}}(\mathbf{x},\mathbf{z},\mathbf{h};s_0, \hat{u}).
\end{align}

\subsection{Cross complementarity and switch detection}\label{sec:fesd_cross_comp}
We proceed with stating the cross complementarity conditions. 
Their objective is to prohibit active set changes within a finite element and to ensure exact switch detection.
\begin{figure}[t]
	\centering
	{\includegraphics[width=0.45\textwidth]{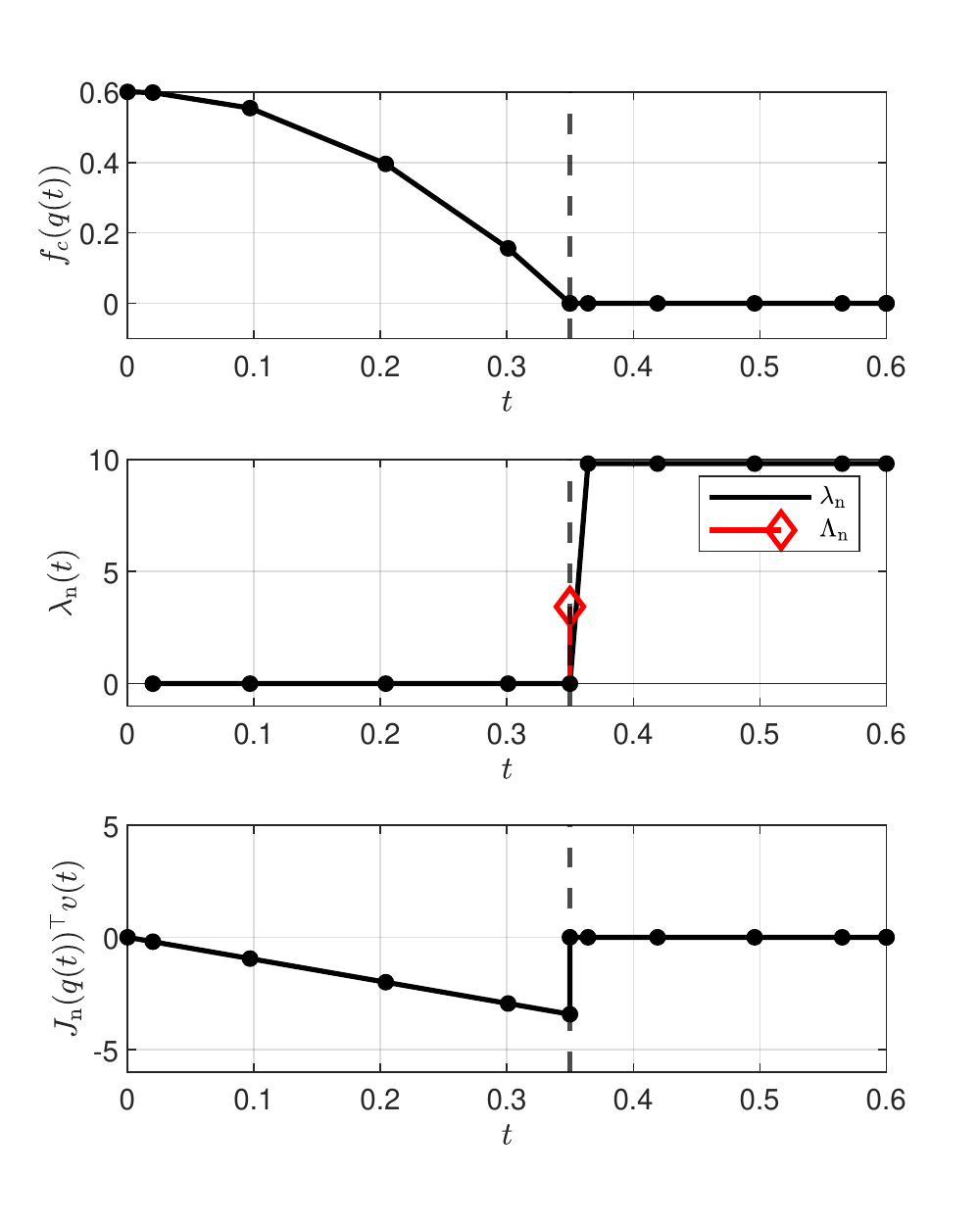}}
	{\includegraphics[width=0.45\textwidth]{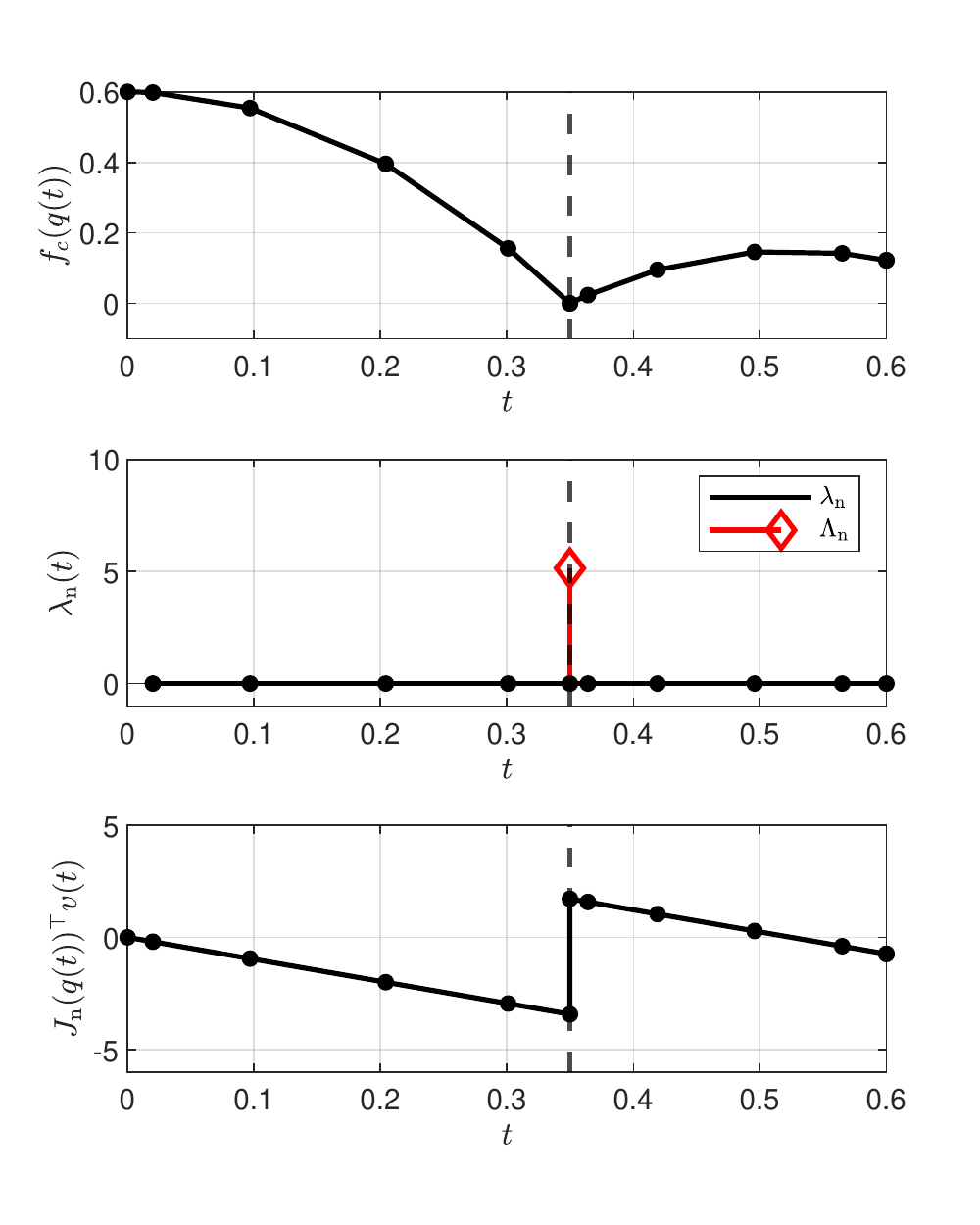}}
	\caption{
		Illustration of trajectories corresponding to an elastic and an inelastic impact.
		The top plots show the position, the middle the normal contact force, and the bottom the normal velocity.
		The left plots correspond to an inelastic impact and the right plots to a partially elastic impact.
		The black markers represent the stage values of a Runge-Kutta Radau IIA method with $\Nstg = 5$.
		The vertical dashed line marks the switching time.
		It can be seen that $\fci(q_{n,j}) > 0$ and $\lambda_{\mathrm{n},n,j'} > 0$  with $j \neq j'$ would lead to the violation of the cross complementarity conditions \eqref{eq:cross_comp_contact}. 
		}
	\label{fig:fesd_cross_comp_contact}
\end{figure}
\subsubsection{Contact making and breaking}
We regard first the discrete-time version of the contact conditions $0\leq \fci(q(t)) \perp \zni(t)\geq 0$.
Recall that under our assumptions, $\fci(q(t))$ is a continuous function of time, whereas  the functions $\zni(t)$ are in general discontinuous.
Consequently, during an active set change in the $i$-th complementarity pair, i.e., during contact breaking or making, we have that $\fci(q(\ts)) = 0$, cf. Figure \ref{fig:fesd_cross_comp_contact}.
These observations motivate, to define in addition to \eqref{eq:irk_ode_alg_contact_comp}, the following \textit{cross complementarity conditions} for CLS
\begin{align}\label{eq:cross_comp_contact}
	0 = \lambda_{\mathrm{n},n,j}^i~\fci (q_{n,j'}),\quad i=1,\ldots,\ncontacts,\ j =1,\ldots, \Nstg, \ j'=0,\ldots, \Nstg.
\end{align}
Note that we also include the boundary point values $\fci (q_{n,0})$, as they play a crucial role in the switch detection, as discussed next.
\begin{lemma}\label{lem:cross_cc_statemnt}
	Regard a fixed $n \in \{0,\ldots,\NFE-1\}$.
	If any $\lambda_{\mathrm{n},n,j}^i$ with $j \in \{1,\ldots, \Nstg\}$ is positive, then all $\fci (q_{n,j'})$ with $j'\in \{0,\ldots, \Nstg\}$ must be zero.
	Conversely, if any $\fci (q_{n,j'})$ is positive, then all $\lambda_{\mathrm{n},n,j}^i$ are zero.
\end{lemma}
\textit{Proof.} Let $\lambda_{\mathrm{n},n,j}^i>0$, and suppose $\fci (q_{n,j}) = 0 $ and $\fci (q_{n,k})>0$ for some $j\in \{1,\ldots, \Nstg\},\ k\in \{0,\ldots, \Nstg\}, k\neq j$, then $\lambda_{\mathrm{n},n,j}^i \fci (q_{n,k})>0$, which contradicts \eqref{eq:cross_comp_contact}, thus all $\fci(q_{n,k})=0,\ k\in \{0,\ldots, \Nstg\}$.
The converse is proven with the same argument. \qed

In the case of an inelastic impact, suppose that some $\fci(q_{n-1,k}) >0$ and that $\lambda_{\mathrm{n},n,j}^i>0$. 
As a consequence of Lemma \ref{lem:cross_cc_statemnt}, it follows that: 
\begin{align*}
	\fci(q_{n,0}) = 0.
\end{align*}
Therefore, we have implicitly a constraint that forces $h_n$ to adapt such that the switch is detected exactly.
The left plots in Figure \ref{fig:fesd_cross_comp_contact} illustrate this switching case.

In the case of elastic impacts, we have before and after the impact a zero contact force, i.e., $\lambda_{\mathrm{n},n-1,j}^i=0$ and $\lambda_{\mathrm{n},n+1,j}^i=0$ for all $j = 1,\ldots, n_q$. 
Consequently, we have $\fci (q_{n-1,j}) \geq 0$ and $\fci (q_{n,j}) \geq 0$. 
However, since $\Znin >0$ in this case, from \eqref{eq:position_switch_detection} it follow that $\fci (q_{n,0}) = 0$ must hold, and the correct $h_n$ is selected.
The right plots in Figure \ref{fig:fesd_cross_comp_contact} illustrate this switching case.

The constraints \eqref{eq:cross_comp_contact} are for the sake of clarity given in their sparsest form.
However, the nonnegativity of  $\lambda_{\mathrm{n},n,j}^i$ and $\fci (q_{n,j})$ allows to aggregate the constraint and obtain equivalent formulations with fewer constraints.

\subsubsection{Switches due to friction}
Making and breaking contacts is not the only type of switching that occurs in a CLS.
If a contact is active, also friction forces act on the body that may lead to further discontinuities.
In particular, during this phase of motion, there may be discontinuous changes in the velocity, i.e., due to transition from slip to stick motion, from stick to slip, or during the slip motions, the velocity may change its direction and thus the friction force changes discontinuously.

These switches will be isolated with another set of cross complementarity constraints.
If no impacts occur, we have that the velocity is continuous, i.e., $v_n = v_{n,0}$.
From \eqref{eq:velocity_switch_detection}, we have that the approximation of Lagrange multiplier $\gamma^i$ corresponding to $v_{n,0}$ is $\Gamma^{i}_n$,
and for ease of notation we define $\gamma^{i}_{n,0} = \Gamma^{i}_n$. 
Recall that $\gamma^i$ is proportional to the magnitude of the tangential velocity, which is a continuous function of time as long as the $i$-th contact stays active.
On the other hand, the friction force $\zti$ might have jump discontinuities if the tangential velocity has a zero crossing.
Therefore, similar to the contact conditions, we have in $0\leq \gamma^i \perp  \beta^i = ((\mu^{i} \lambda_{\mathrm{n}}^{i})^2 - \| \lambda_{\mathrm{t}}^{i} \|_2^2) \geq 0$ a complementarity between a continuous and discontinuous function of time.
In order to prohibit active set-changes within the finite elements, similar to \eqref{eq:cross_comp_contact}, we introduce the following set of cross complementarity constraints for $i=1,\ldots,\ncontacts$:
\begin{align}\label{eq:cross_complementarity_friction_1}
	0 = \gamma_{n,j'}^{i}~ \beta_{n,j}^{i} , \quad j =1,\ldots, \Nstg, \ j'=0,\ldots, \Nstg.
\end{align}
For the approximation of the continuous-time variable $\gamma^i$, we include the boundary point value $\gamma^{i}_{n,0}$.
If we have some $\gamma^{i}_{n-1,j} >0$ and $\beta_{n,j}^{i} = (\mu^{i} \lambda_{\mathrm{n},n,j}^{i})^2 - \| \lambda_{\mathrm{t},n,j}^{i} \|_2^2 = 0$, then by the same arguments as in Lemma~\ref{lem:cross_cc_statemnt}, we have that
$\gamma_{n,0}^{i} =  \Gamma_{n}^{i} = 0$.
From \eqref{eq:velocity_switch_detection} it follows that
\begin{align*}
      \Jti(q_n)^\top v_{n,0}  = 0,
\end{align*}
which implicitly provides a condition that isolates the zero crossing of the tangential velocity.
\begin{figure}[t]
	\centering
	{\includegraphics[width=0.32\textwidth]{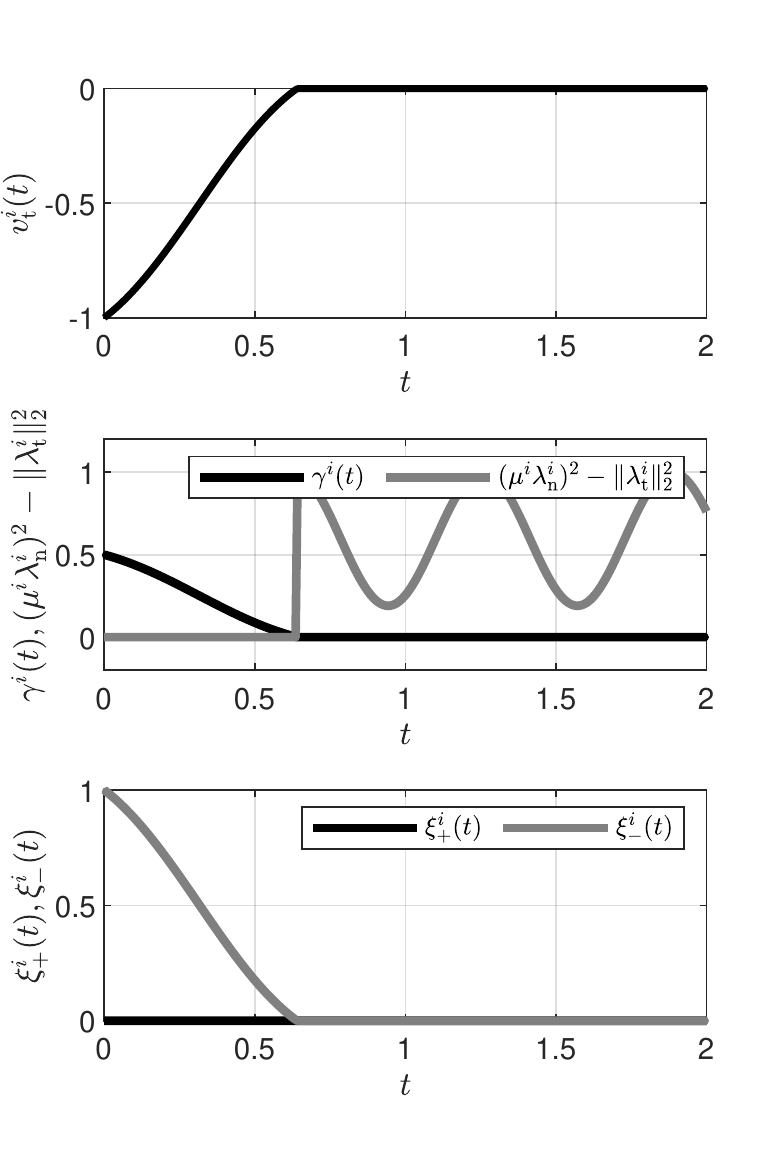}}
	{\includegraphics[width=0.32\textwidth]{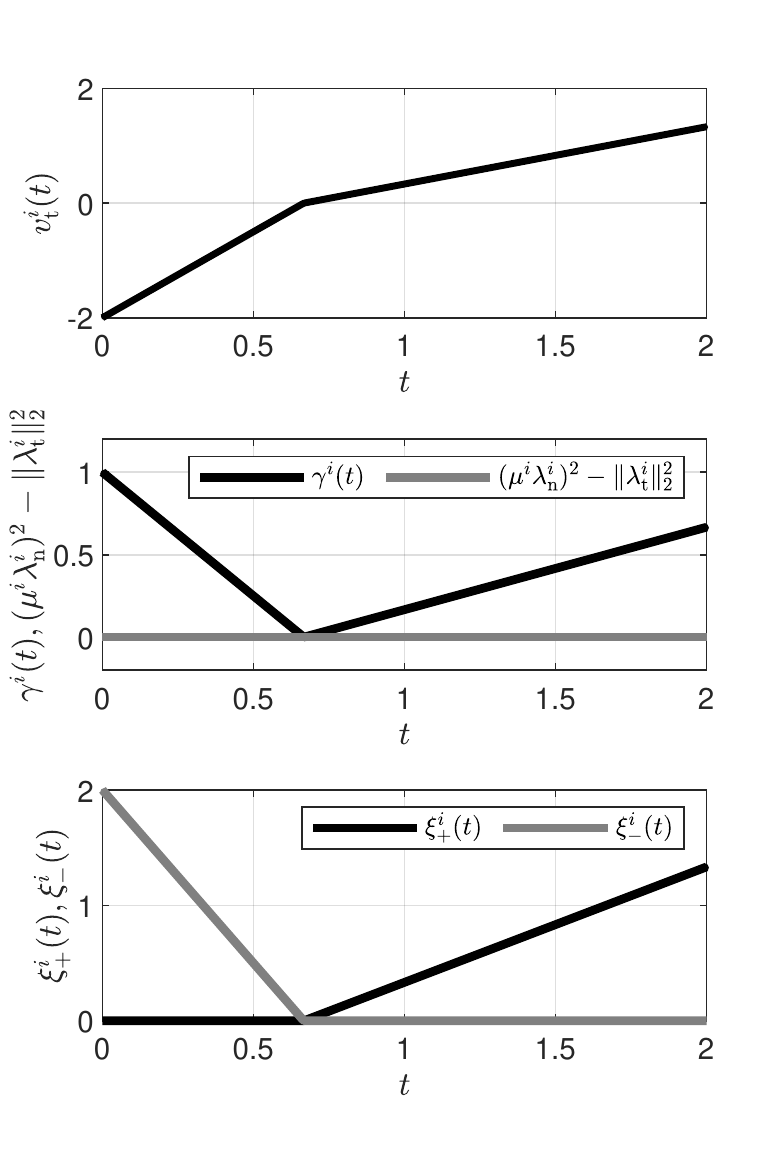}}
	{\includegraphics[width=0.32\textwidth]{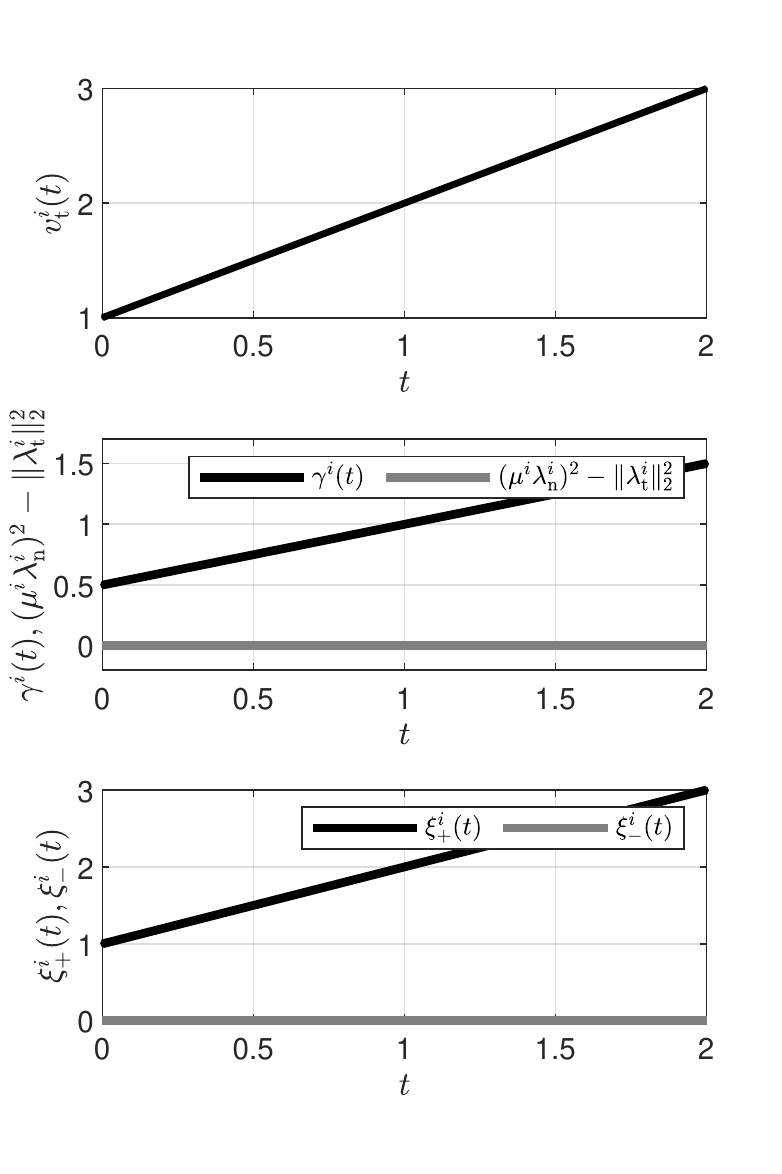}}
	\caption{Different switching cases due to friction. 
		The left plot shows a particle that transitions from slip to stick motion, where the friction force during the stick phase is variable.
		The middle plot shows a discontinuity where a particle remains in slip motion.
		The right plot shows a particle in slip motion with no switches.
	}
	\label{fig:fesd_friction}
\end{figure}

These conditions are sufficient to detect the switches in the cases of a transition from slip to stick mode or from stick to slip, cf. Figure \ref{fig:fesd_friction} (left column).
However, they are not sufficient for detecting a zero crossing if the friction force changes discontinuously while the body stays in a slipping motion.
This situation is illustrated in the middle column of Figure~\ref{fig:fesd_friction}, where a particle with a negative initial velocity in slipping motion.
Moreover, a force larger in magnitude than the friction force acts in the opposite direction of the particle's motion.
The particle will at some point reach zero velocity and continue to slip in the opposite direction, now with a smaller acceleration, since the external force and friction act in opposite directions.
Compared to smooth slip motion with no switches (cf. middle plots in second and third column of Figure~\ref{fig:fesd_friction}), the active sets in $0 \leq \gamma^{i}_{n,j} \perp \beta^{i}_{n,j} \geq 0, j = 1,\ldots,\Nstg$ from \eqref{eq:irk_ode_alg_friction1} remain unchanged.
Consequently, the conditions in \eqref{eq:cross_complementarity_friction_1} are trivially satisfied and cannot isolate the switch.
To over come this problem, one can use the positive and negative parts ($\xi_+^i$ and $\xi_-^i$) of the tangential velocities with Eq. \eqref{eq:friction_equations_3D_augmented}.
We propose the following set of cross complementarity conditions:
\begin{align}\label{eq:cross_complementarity_friction_2}
	0 = \mathrm{diag}(\xi^{i}_{+,n,j'}) \xi^{i}_{-,n,j}, \quad j =0,\ldots, \Nstg, \ j'=0,\ldots, \Nstg.
\end{align}
It can be seen from the bottom plots in the second and third column of Figure~\ref{fig:fesd_friction}, that the active sets in \eqref{eq:friction_pos_neg_comp} are not the same.
Therefore, the cross complementarity conditions \eqref{eq:cross_complementarity_friction_2} are not implied by the standard stage-wise complementarity conditions \eqref{eq:friction_pos_neg_discrete}.
By similar reasoning, we have in the described switching case that $\xi^{i}_{+,n,j'} = \xi^{i}_{-,n,j} = 0$, and from \eqref{eq:friction_pos_neg_def} we obtain that:
\begin{align*}
	\Jti(q_n)^\top v_{n,0}  = 0,
\end{align*}
which makes sure that the switches are detected also if the body remains in slipping motion.

We collect all cross complementarity conditions \eqref{eq:cross_comp_contact}, \eqref{eq:cross_complementarity_friction_1}  and \eqref{eq:cross_complementarity_friction_2}	into the equation:
\begin{align}\label{eq:cross_comp}
	0 = G_{\mathrm{cross}}(\mathbf{x},\mathbf{z};s_0).
\end{align}

\subsection{Step equilibration}
If two neighboring finite elements have the same active set in the corresponding complementarity conditions, then the cross complementarity conditions are implied by the standard stage-wise complementarity conditions~\cite{Nurkanovic2022}.
Consequently, we obtain spurious degrees of freedom in the step sizes $h_n$.
We overcome this by introducing additional constraints that appear only in this case, and that remove these degrees of freedom.

For this purpose we need to introduce an indicator function $\eta_n(\mathbf{x},\mathbf{z})$, which is strictly positive if no switch happens at the $n$-th grid point and zero if a switch occurs.
The step equilibration equations read as
\begin{align}\label{eq:step_eq}
	0&= G_{\mathrm{eq}}(\mathbf{x},\mathbf{z},\mathbf{h};s_0)
	\coloneqq 
	\begin{bmatrix}
		(h_{1}-h_{0})\eta_1(\mathbf{x},\mathbf{z}) \\
		\vdots\\
		(h_{\NFE\!-\!1}-h_{\NFE\!-\!2})\eta_{\NFE-1}(\mathbf{x},\mathbf{z})
	\end{bmatrix}.
\end{align}
If no switches happen, these conditions imply an equidistant grid. 
If switches occur at some grid points $t_n$, then the resulting grid will be piecewise equidistant.

In the remainder of this section, we define the functions $\eta_n$, which encode the switching logic.
We start with the switches that happen by making and breaking contact.
Following \cite{Nurkanovic2022}, we define for $n = 1,\ldots, \NFE$ forward and backward sums for the normal contact force and gap function for every $i = 1,\ldots,\ncontacts$:
\begin{subequations}
\begin{align}
	&\sigma^{\fci,\mathrm{B}}_{n} = \sum_{j=0}^{\Nstg} \fci(q_{n-1,j}), \quad
	\sigma^{\fci,\mathrm{F}}_{n} = \sum_{j=0}^{\Nstg} \fci(q_{n,j}),
	\\
	&\sigma^{\zni,\mathrm{B}}_{n} = \sum_{j=1}^{\Nstg} \lambda^{i}_{\mathrm{n},n-1,j}, \quad
	\sigma^{\zni,\mathrm{F}}_{n} = \sum_{j=1}^{\Nstg} \lambda^{i}_{\mathrm{n},n,j}.
\end{align}
\end{subequations}
The values of these sums change qualitatively before and after a switch, which enables us to encode the switching logic.
For example, in the case of making contact at $t_n$, in the inelastic case we have $\sigma^{\fci,\mathrm{B}}_{n}> 0$, $ \sigma^{\fci,\mathrm{F}}_{n} = 0$ and
$\sigma^{\zni,\mathrm{B}}_{n} = 0, \sigma^{\zni,\mathrm{F}}_{n} > 0$. 
In the case of elastic impacts, only $\fci(q_{n,0}) = 0$ and $\sigma^{\fci,\mathrm{B}}_{n}>0$, $\sigma^{\fci,\mathrm{F}}_{n} >0$, cf. Figure~\ref{fig:fesd_cross_comp_contact}.
With logical \textit{or} and \textit{and} operations on these quantities we can indicate whether there was a switch or not.
The contact conditions permit four different switching scenarios, which are listed in Table \ref{tab:switching_logic_1}.
For ease of readability, positive values of a variable are represented by one, and zero values are represented by zero.
\begin{table}[t]
	\centering
	\caption{Switching logic for contact making and breaking in the case of elastic impacts. 
		One indicates that a variable is strictly positive and a zero that it is zero.}
		\begin{tabular}{|c|ccc|cc|cc|c|c|}
			\hline
			\textbf{Switching case} & $\fci(q_{n,j})$ & $\sigma^{\fci,\mathrm{B}}_{n}$  & $\sigma^{\fci,\mathrm{F}}_{n}$ &  $\sigma^{\zni,\mathrm{B}}_{n}$  & $\sigma^{\zni,\mathrm{F}}_{n}$ & $\pi_n^{\fci}$ & $\pi_n^{\zni}$ &$\kappa_n^{i}$&Switch at $t_n$\\
			\hline
			Contact making     &0 &1 &1  &0 &0   &0 &0   &0 &yes\\ 
			Contact breaking   &0 &0 &1  &1 &0   &0 &0   &0 &yes\\ 
			Free flight        &1 &1 &1  &0 &0   &1 &0   &1 &no\\ 
			Persistent contact &0 &0 &0  &1 &1   &0 &1   &1 &no\\ 
			\hline
		\end{tabular}
		\label{tab:switching_logic_1}
	\end{table}
To mimic logical \textit{and} operations for the gap functions we define:
\begin{align*}
	\pi_n^{\fci}& = \fci(q_{n,0})\sigma^{\fci,\mathrm{B}}_{n}\sigma^{\fci,\mathrm{B}}_{n}.
\end{align*}
If there is persistent contact, all terms in the product are zero, and it follows that $\pi_n^{\fci} = 0$.
Similarly, we define the product of the sums of contact forces:
\begin{align*}
	\pi_n^{\zni}& = \sigma^{\zni,\mathrm{B}}_{n}\sigma^{\zni,\mathrm{B}}_{n}.
\end{align*}
In the persistent contact case, $\pi_n^{\zni}$ is strictly positive.
In case of free flight, we have the symmetric case of $\pi_n^{\fci} >0 $ and $\pi_n^{\zni} = 0$.
If neither contact breaking or making happens, one of the two expressions is zero and one is positive.
On the other hand, if we have contact making or breaking, both products are zero.
This can be summarized with a logical \textit{or} operation, which is modeled with 
\begin{align}
	\kappa_n^i =  \pi_n^{\fci}+\pi_n^{\zni}.
\end{align}
If no switch occurs $\kappa_n^i$ is positive, and otherwise it is zero.
All switching cases are enumerated in Table~\ref{tab:switching_logic_1}.
For inelastic impacts, the reasoning is only different for the case of contact making.
In this case, we have that $\fci(q_{n,j})= 0$, $\sigma^{\fci,\mathrm{B}}_{n} >0$ and $\sigma^{\fci,\mathrm{F}}_{n} =0$. F
For the backward and forward sums for the contact forces it holds that $\sigma^{\zni,\mathrm{B}}_{n} = 0$ and $\sigma^{\zni,\mathrm{B}}_{n}>0$. 
Thus, we obtain that $\pi_n^{\fci} = 0, \pi_n^{\zni} = 0$ and $\kappa_n^i =  0$.

Next, we define similar functions for switches due to the nonsmooth friction forces.
We have qualitatively five different scenarios, three with a switch and two without.
They are listed in Table \ref{tab:switching_logic_friction}. 
We assume in the table for the slip-stick transition positive tangential velocity before the switch, and for the stick to slip negative tangential velocity after the switch. 
All other combinations follow the same reasoning and we omit stating them explicitly for brevity.

To encode the switching logic, we follow a similar approach and define the following forward and backward sums:
	\begin{align*}
		&\sigma^{\mathrm{\xi_+^i,B}}_{n,k} = \sum_{j=0}^{\Nstg} \xi^i_{+,n-1,j},\quad
		\sigma^{\xi_{+}^i,\mathrm{F}}_{n,k} = \sum_{j=0}^{\Nstg} \xi^i_{+,n,j}, &k = 1,\ldots,n_{\mathrm{t}},\\
		&\sigma^{\mathrm{\xi_-^i,B}}_{n,k} = \sum_{j=0}^{\Nstg} \xi^i_{-,n-1,j}, \quad
		\sigma^{\xi_{-}^i,\mathrm{F}}_{n,k} = \sum_{j=0}^{\Nstg} \xi^i_{-,n,j},&k = 1,\ldots,n_{\mathrm{t}},\\	
		&\sigma^{\beta^i,\mathrm{B}}_{n} = \sum_{j=1}^{\Nstg} \beta^i_{n-1,j}, \quad
		\sigma^{\beta^i,\mathrm{F}}_{n,k} = \sum_{j=1}^{\Nstg} \beta^i_{n,j}.
	\end{align*}
We omitted to define forward and backwards sums for the variables $\gamma^i_{n,j}$, as the three above are sufficient to cover all switching cases.
For the positive and negative tangential velocity parts, we have for every component a separate sum. 
In the planar case, it is simply $n_{\mathrm{t}} = 1$ and in the three-dimensional case we have $n_{\mathrm{t}} = 2$. 
For brevity, Table \ref{tab:switching_logic_friction} reports the cases with $n_{\mathrm{t}} = 1$, but the reasoning is similar for $n_{\mathrm{t}} = 2$.

Similar to contact making and breaking, we define variables that relate the sums before and after a grid point $t_n$
\begin{subequations}\label{eq:logical_friction}
\begin{align}
	&\pi_{n,k}^{\xi^i_+} = \sigma^{\xi_{+}^i,\mathrm{B}}_{n,k} \sigma^{\xi_{+}^i,\mathrm{F}}_{n,k}, &k = 1,\ldots,n_{\mathrm{t}},\\
	&\pi_{n,k}^{\xi^i_-} = \sigma^{\xi_{-}^i,\mathrm{B}}_{n,k} \sigma^{\xi_{-}^i,\mathrm{F}}_{n,k}, &k = 1,\ldots,n_{\mathrm{t}},\\
	&\pi_n^{\beta^i} = \sigma^{\beta^i,\mathrm{B}}_{n} \sigma^{\beta^i,\mathrm{F}}_{n}.
\end{align}
\end{subequations}
For sake of illustration, we discuss the example of a transition from slip to stick mode in the two dimensional case. 
Suppose that the tangential velocity before the switch is positive. 
In this case, we have that 
$\sigma^{\xi_{+}^i,\mathrm{B}}_{n} > 0$, $\sigma^{\xi_{+}^i,\mathrm{F}}_{n} = 0$ and 
$\sigma^{\xi_{-}^i,\mathrm{B}}_{n} = 0, \sigma^{\xi_{-}^i,\mathrm{F}}_{n} = 0$.
During slip motion, the friction force had its maximum value, hence $\sigma^{\beta^i,\mathrm{B}}_{n} = 0$.
After the switch during stick motion motion, the friction force might be at its maximum or not, hence $\sigma^{\beta^i,\mathrm{B}}_{n} \geq 0$.
Therefore, all variables defined in \eqref{eq:logical_friction} are zero.

To summarize, the effects of all switches in the frictional case, we define
\begin{align*}
	\zeta_n^i = (\sigma^{\fci,\mathrm{B}}_{n}+\sigma^{\fci,\mathrm{F}}_{n})+ (\sum_{k=1}^{n_{\mathrm{t}}} (\pi_{n,k}^{\xi^i_+} + \pi_{n,k}^{\xi^i_-}) + \pi_n^{\beta^i}).
\end{align*}
The first term makes sure that $\zeta_n^i$ can be zero only if the $i$-th contact is active before and after $t_n$.
Only in this case the first term is zero, otherwise it is strictly positive,  cf. Table \ref{tab:switching_logic_1}. 
This avoids that the step equilibration for the frictional switches has any influences for inactive contacts.
If any of the terms in the second parentheses is positive, then there was no switch, i.e., $\zeta_n^i>0$. 
The sign and switching logic for all cases is summarized in Table~\ref{tab:switching_logic_friction}.
The validity of the entries in Table~\ref{tab:switching_logic_friction} can be easily verified by looking at Figure~\ref{fig:fesd_friction}.

It is left to summarize the effects of all contacts.
For a single contact, we introduce the switching logic via
\begin{align*}
	\nu^i_n = \kappa_n^i\,\zeta_n^i.
\end{align*}
The variable $\nu^i_n$ is zero if there is a switch due to contact making or breaking or due to friction.
If there are no switches it is positive.
The step equilibration should not be present if there is a switch in any of the contact. 
This is modeled by defining:
\begin{align}
	\eta_n^i = \prod_{i = 1}^{\ncontacts} \nu^i_n.
\end{align}
This completes the definition of the step equilibration conditions~\eqref{eq:step_eq}.
\begin{table}[t]
	\centering
	\caption{Switching logic in contact making and breaking in the case of elastic impacts. In the table there is a one if a variable is strictly positive and a zero if it is zero.
		It is assumed that $\sigma^{\fci,\mathrm{B}}_{n}+\sigma^{\fci,\mathrm{F}}_{n} = 0$.
	}
	\begin{tabular}{|c|cc|cc|cc|ccc|c|c|}
		\hline
		\textbf{Switching case} 
		&$\sigma^{\xi_{+}^i,\mathrm{B}}_{n}$ &$\sigma^{\xi_{+}^i,\mathrm{F}}_{n}$
		&$\sigma^{\xi_{-}^i,\mathrm{B}}_{n}$ &$\sigma^{\xi_{-}^i,\mathrm{F}}_{n}$
		&$\sigma^{\beta^i,\mathrm{B}}_{n}$ &$\sigma^{\beta^i,\mathrm{F}}_{n}$
		& $\pi_n^{\xi^i_+}$
		& $\pi_n^{\xi^i_-}$
		& $\pi_n^{\beta^i}$ 
		&$\zeta_n^i$& Switch \\
		\hline
		Switch during slip   &1 &0  &0 &1  &0 &0  &0&0&0  &0& yes \\ 
		Slip to stick    	 &1 &0  &0 &0  &0 &1  &0&0&0  &0 &yes \\ 
		Stick to slip   	 &0 &0  &0 &1  &1 &0  &0&0&0  &0 &yes \\ 
		Stick phase   	     &0 &0  &0 &0  &1 &1  &0&1&1  &1 &no \\ 
		Slip motion 	 	 &1 &1  &0 &0  &0 &0  &1&0&0  &1 &no \\ 
		\hline
	\end{tabular}
	\label{tab:switching_logic_friction}
\end{table}
\begin{remark}(Heuristic step equilibration)
The term $\eta_n^i$ involves products of many different terms and it may take values in a broad range. 
For a better scaling we may use $\tanh(\eta_n^i/ \epsilon_{\eta})$, with $\epsilon_{\eta}\in (0, 1]$.
Moreover, since $\eta_n^i$  might be quite nonlinear, we propose to use the following heuristic approach.
Step equilibration can approximately be achieved by adding the term $\rho_h \sum_{n = 0}^{\NFE-2}(h_{n+1}-h_n)^2$ to the cost function, where $\rho_h > 0$ is a weighting factor. 
In optimal control problems, one should not choose a too large value for $\rho_h$ relative to the other objective terms, as this might introduce a bias towards controls that result in an equidistant discretization grid.
\end{remark}

\subsection{FESD-J - summary}
We have now introduced all equations that define the FESD-J discretization of the CLS \eqref{eq:cls}. 
Given an initial state $s_0$, these equations provide a numerical approximation of the differential and algebraic states over the time interval $[0,T]$.
We collect the FESD equations \eqref{eq:rk_equations}, \eqref{eq:cross_comp} \eqref{eq:step_eq}, and \eqref{eq:sum_h} in the discrete-time system
\begin{subequations} \label{eq:fesd_compact}
	\begin{align}
		s_{1} &=F_{\fesd}(\mathbf{x}), \label{eq:fesd_compact_state_transition_step_representation}\\
		0&= G_{\fesd}(\mathbf{x},\mathbf{z},\mathbf{h};s_0, \hat{u}, T)
		\coloneqq
		\begin{bmatrix}
			G_{\mathrm{rk}}(\mathbf{x},\mathbf{z},\mathbf{h};s_0, \hat{u})\\
			G_{\mathrm{cross}}(\mathbf{x},\mathbf{z};s_0)\\
			G_{\mathrm{eq}}(\mathbf{x},\mathbf{z},\mathbf{h};s_0)\\
			\sum_{n=0}^{\NFE-1} h_n - T
		\end{bmatrix},
	\end{align}
\end{subequations}
where the state transition map $F_{\fesd}(\mathbf{x}) = (q_{\NFE},v_{\NFE})$ provides the state approximation $s_1 \approx x(T)$, and $G_{\fesd}(\mathbf{x},\mathbf{z},\mathbf{Z},\mathbf{h};s_0, \hat{u}, T)$ collects all other internal computations.
Here, the control variable $\hat{u}$, the interval length $T$, and the initial value $s_0$ are given parameters.

\begin{remark}(Switching time interpretation)
An alternative interpretation for the variable step sizes $h_n$ is that they are discrete time control variables that control the interval lengths and thus implicitly determine the switching times.
This can be seen by interpreting $t$ as a clock state, with the dynamics 
\begin{align*}
	&\dot{t} = h,\ t(0) = 0,\ t(1) = T,
\end{align*}
whose time discretization reads as:
	\begin{align*}
		&t_{n+1} = t_n  + h_n,\quad n = 0,\ldots,\NFE-1,\\ 	
		&t_0 = 0,\	t_{\NFE} = T.
	\end{align*}
The control variables $h_n$ are determined implicitly via cross complementarity and step equilibration constraints.
\end{remark}

\section{Direct optimal control}\label{sec:fesd_ocp}
Our main motivation for introducing \FESDj\ is to discretize and solve Optimal Control Problems (OCPs) subject to CLS.
A continuous time OCP reads as 
\begin{subequations} \label{eq:ocp}
	\begin{align}
		\min_{x(\cdot),\lambda(\cdot), u(\cdot)} \quad & \int_{0}^{T} L(x(t),u(t))\dd t +  L_{\mathrm{t}}(x(T)) \\
		\textrm{s.t.} \quad  &x(0) = {x}_0, \\
		&\textrm {Eq}~\eqref{eq:cls},\ \textrm{a.a. } t \in [0, T], \label{eq:ocp_cls} \\
		&0\geq G_{\mathrm{p}}(x(t),u(t)),\ t \in [0,T],\\
		&0\geq G_{\mathrm{t}}(x(T)),
	\end{align}
\end{subequations}
where ${x}_0$ is a given initial value, $L: \R^{n_x} \times \R^{n_u} \to \R$ is the running cost and $L_{\mathrm{t}}:\R^{n_x}\to \R$ is the terminal cost, $s_0\in\R^{n_x}$ is a given initial value.
The path and terminal constraints are defined by the functions $G_{\mathrm{p}} : \R^{n_x}  \times \R^{n_u} \to \R^{n_{\mathrm{p}}}$ and $G_{\mathrm{t}} : \R^{n_x}  \to \R^{n_{\mathrm{t}}}$, respectively.

Now we consider a direct transcription of \eqref{eq:ocp}, with the \FESDj\ method from the previous section.
In particular, we make use of the short discrete-time system notation in \eqref{eq:fesd_compact} for the time-discretization of \eqref{eq:ocp_cls}.
Consider $\Nctrl\geq 1$ control intervals of equal length, which are indexed by $k$.
We take a piecewise constant control discretization, where the control variables are collected $\mathbf{u} = (\hat{u}_0,\ldots,\hat{u}_{\Nctrl-1})\in \R^{\Nctrl n_u}$.
On each control interval $k$, we use a \FESDj\ discretization with $N_{\mathrm{FE}}$ internal finite elements.
The state values at the control interval boundaries are collected in the vector  $\mathbf{s} = (s_0,\ldots,s_{\Nctrl})\in\R^{(\Nctrl+1)n_x}$.
We collect in ${\mathcal{Z}} = (\mathbf{z}_0,\ldots,\mathbf{z}_{\Nctrl-1})$, ${\mathcal{X}} = (\mathbf{x}_0,\ldots,\mathbf{x}_{\Nctrl-1})$ all internal variables, and in $\mathcal{H} = (\mathbf{h}_0,\ldots,\mathbf{h}_{\Nctrl-1})$ all step sizes.

The discrete-time variant of \eqref{eq:ocp} is given as
\begin{subequations}\label{eq:ocp_discrete_time}
	\begin{align}
		\min_{\mathbf{s},\mathbf{u},\mathcal{Z},\mathcal{X},\mathcal{H}} \quad & \sum_{k=0}^{\Nctrl-1} \hat{L}(s_k,\mathbf{x}_k,q_k)+ L_{\mathrm{t}}(s_{\Nctrl}) \\
		\textrm{s.t.} \quad  &s_{0} = {x}_0,\\
		&{s}_{k+1}  = F_{\fesd}(\mathbf{x}_k),  &k = 0,\ldots,\Nctrl-1,\\	
		&0 = G_{\fesd}(\mathbf{x}_k,\mathbf{z}_k,\mathbf{h}_k;s_k, \hat{u}_k, \frac{\Tctrl}{\Nctrl}),  &k = 0,\ldots,\Nctrl-1,\label{eq:ocp_discrete_time_cc}\\
		&0 \geq G_{\mathrm{p}}(s_k,\hat{u}_k),  &k = 0,\ldots,\Nctrl-1,\\	
		&0 \geq G_{\mathrm{t}}(s_{\Nctrl}),
	\end{align}
\end{subequations}
where $\hat{L}:\R^{n_x}\times \R^{(\NFE+1)\Nstg n_x} \times \R^{n_u}\to \R$ is the discretized running costs.
Since \eqref{eq:ocp_discrete_time_cc} contains all complementarity constraints appearing in the \FESDj\ discretization, the nonlinear program \eqref{eq:ocp_discrete_time} is a mathematical program with complementarity constraints.
In NOSNOC \cite{Nurkanovic2022b}, MPCCs are solved by solving a sequence of related and relaxed NLPs within a homotopy approach~\cite{Anitescu2007,Scholtes2001}.

\section{Illustrative numerical examples}\label{sec:examples}
In this section, we demonstrate the use of the FESD-J method on simulation and optimal control examples.
All presented examples are publicly available in NOSNOC v0.4.0\footnote{\url{https://github.com/nurkanovic/nosnoc/releases/tag/v0.4.0}}.
\FloatBarrier
\subsection{Two balls connected by a spring -- an integration order experiment}
\begin{figure}[]
	\centering
	\includegraphics[width=.99\textwidth]{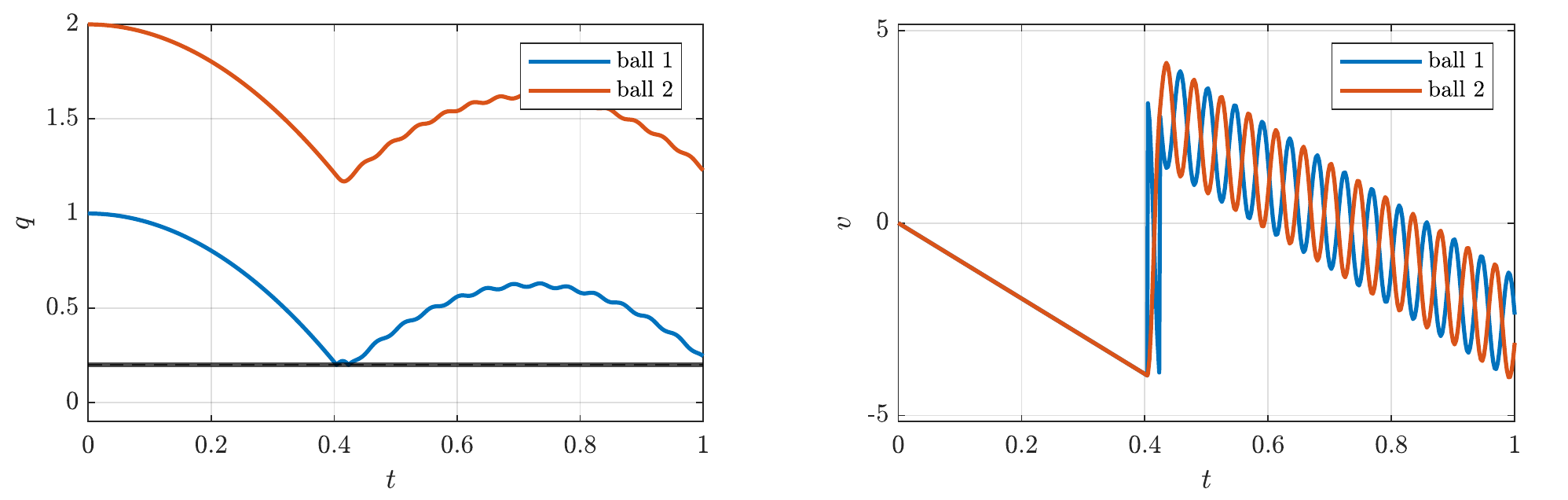}
	\caption{Trajectory of two balls connected by a spring.}
	\label{fig:two_balls_traj}
\end{figure}

In this section, we demonstrate the integration order of the proposed discretization scheme experimentally, on an example of two balls connected by a spring, which has been studied before in \cite{Chen2013}.
The model consists of $q = [q_1, q_2], v=[v_1, v_2] \in \R^2 $, where $q_i, v_i$ denote the position and velocity of ball $i\in\{1, 2\} $, respectively.
The radii of both balls is $R = 0.2$ and the first ball can make contact with the ground, which is modeled with the gap function
\begin{align}
	f_c(q) = q_1  - R \geq 0.
\end{align}
The impact is partially elastic with $\CoR = 0.8$.
The gravitational acceleration is $ g = 9.81 $, the stiffness of the spring connecting the balls is $ k = 10^4$ and its rest length is $l=1$.
The mass matrix is constant $M = \mathrm{diag}(m_1, m_2)$, where $m_1= 1$ and $m_2 =1 $ are the masses of the two balls.
The CLS reads
\begin{align}
	&\dot{q} = v,\\
	&M \dot{v} = \begin{bmatrix} -m_1 g+ k (q_2-q_1-l) \\ -m_2 g-k (q_2-q_1-l) \end{bmatrix} + J_\mathrm{n}(q)\zn,\\
	&0 \leq \zn  \perp \fc(q) \geq 0.
\end{align}
The system is simulated with initial zero velocities and initial position $q_{1}(0) = 1, q_{2}(0) = 2 $ for a time interval of $ T = 1.0$.
Figure~\ref{fig:two_balls_traj} shows the resulting trajectories.

\begin{figure}[t]
	\centering
	\includegraphics[width=.49\textwidth]{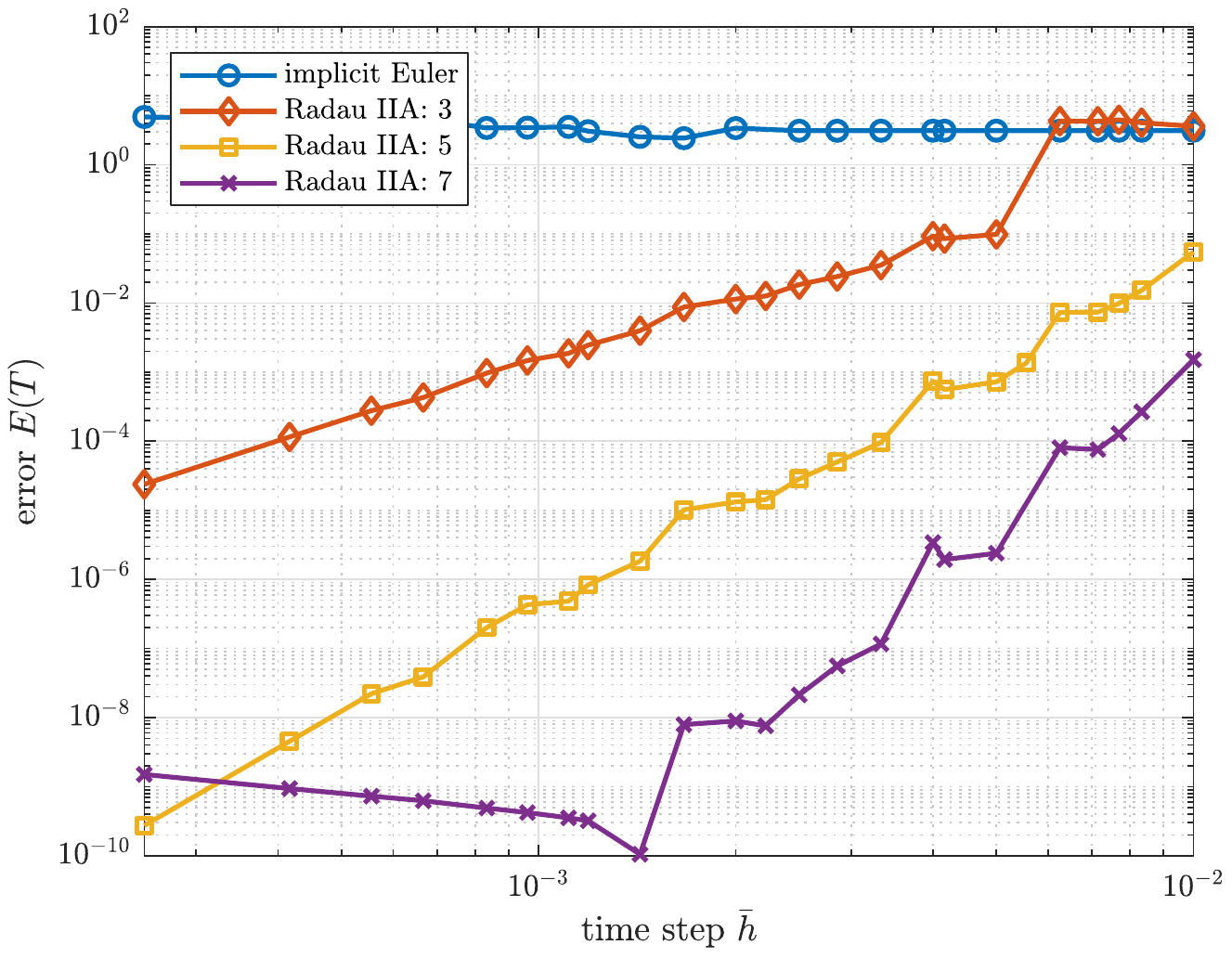}
	\includegraphics[width=.49\textwidth]{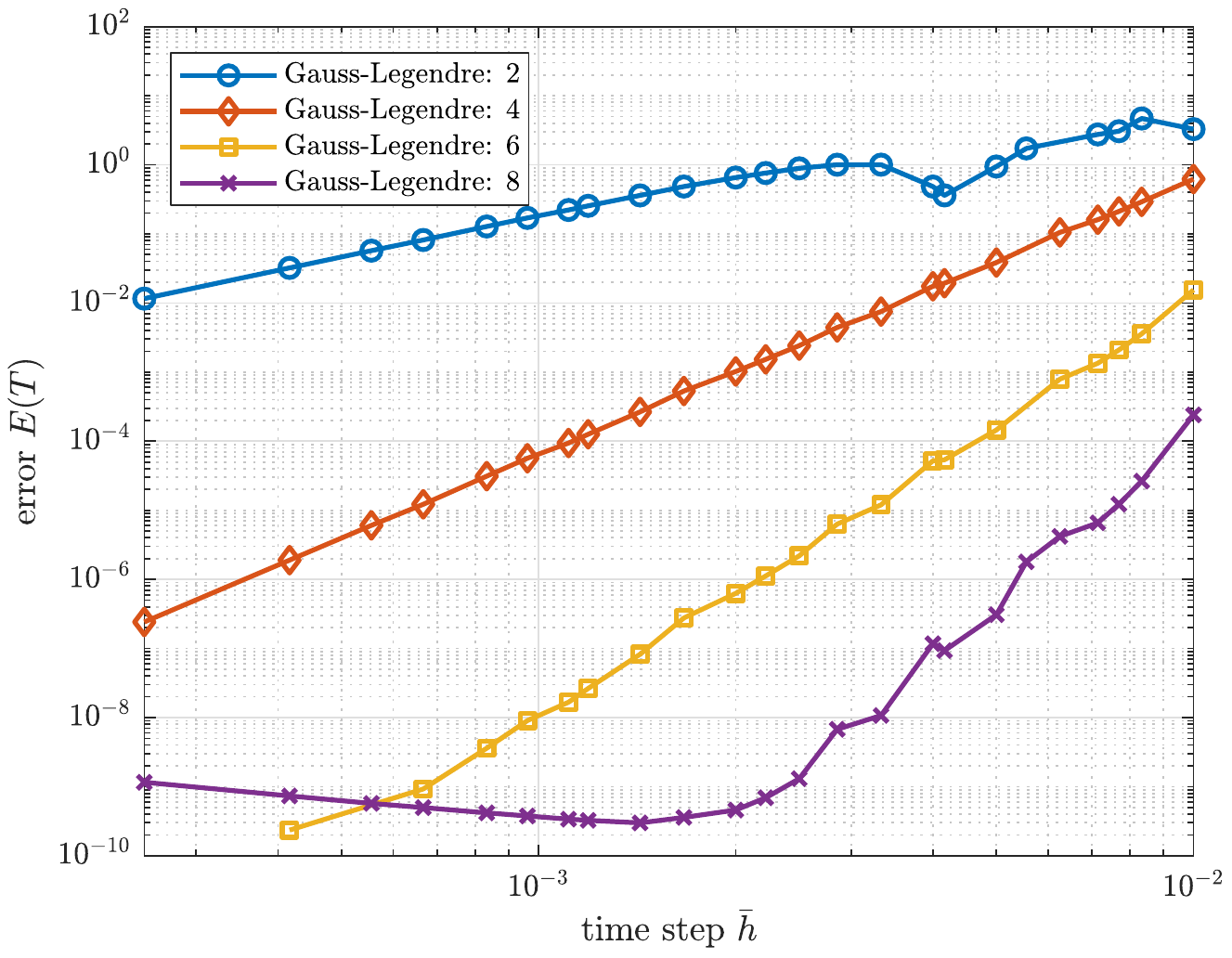}
	\caption{Experimental integration order simulating the scenario depicted in Figure~\ref{fig:two_balls_traj}.}
	\label{fig:order_plot_two_balls}
\end{figure}

In the simulation experiment, we vary the number of simulation steps $N_{\mathrm{sim}} $ and fix the number of finite Elements to $\NFE=2$.
We use the proposed FESD-J method with the Gauss-Legendre and Radau IIA  Butcher tableaus with $\Nstg\in \{1,2,3,4\}$ to simulate the system.
In Figure~\ref{fig:order_plot_two_balls}, we plot the solution accuracy against the time step of one simulation step $\bar{h} = \frac{T}{N_{\mathrm{sim}}  \NFE }$.
Since solving this simulation problem numerically with the existing homotopy methods in NOSNOC can be challenging, in case of indefeasibly, we attempt to solve the subproblem with a modified initial guess, with the impact impulse $\Lambda_\mathrm{n} > 0$.
This improves the overall solutiont robustness significantly.
Moreover, we exclude runs, where any of the $N_{\mathrm{sim}}$ subproblems did not converge to the required accuracy, which make up $ 4.17 \%$ and $ 5.2\%$ of runs for Radau IIA and Gauss-Legendre, respectively.
The plot shows that the established order properties of implicit RK methods are preserved when using them within \FESDj.

\FloatBarrier
\subsection{Manipulation task - inelastic impacts}
In the first optimal control example, we regard two discs that lie in a two-dimensional plane. 
Only one disc can be controlled by a thrust force, whereas the second can only be moved by making inelastic contact (with $\CoR = 0$) with the first one.
The goal is that the discs swap their position with minimal control effort. 
This can be modeled with the following optimal control problem:
\begin{align*}
	\min_{x(\cdot),\zn(\cdot),u(\cdot),} \quad &  \int_{0}^{\Tctrl}\| x(t)-x_\mathrm{r}\|_Q^2+\|u(t)\|_R^2\;\dd t
	+\| x(\Tctrl)-x_\mathrm{r}\|_{Q_T}^2 \\
	\textrm{s.t.} \quad  &x(0) = {x}_0,\\
	&\dot{q}(t) = v(t), \; &t\in [0,\Tctrl], \\
	&M\dot{v}(t) = \begin{bmatrix}
		u(t) - c_{\mathrm{f}} \frac{v_1}{\|v_1+\epsilon_{\mathrm{f}}\|} \\ 
		\mathbf{0}_{2,1} - c_{\mathrm{f}} \frac{v_2}{\|v2+\CoR{\mathrm{f}}\|} 
	\end{bmatrix}+ J_{\mathrm{n}}(q(t)) \zn(t),\ &t\in [0,\Tctrl], \\
	&0 \leq \zn(t) \perp \fc(q)\geq0, \; &t\in [0,\Tctrl], \\
	&0 = J_{\mathrm{n}}(q(t))^\top v(\ts^+),\ \textrm{if} \; f_c(q(\ts)) = 0 \textrm{ and } n(q(t))^\top v(\ts^-)<0,\\
	&-10e \leq q(t) \leq 10e, \; &t\in [0,\Tctrl], \\
	&-5e  \leq v(t) \leq 5e, \; &t\in [0,\Tctrl], \\
	&-30e \leq u(t) \leq 30e, \; &t\in [0,\Tctrl], \\
	&0 \leq \|q_1(t)\|^2-(r_{\mathrm{ob}}+r_1)^2, \; &t\in [0,\Tctrl], \\
	&0 \leq \|q_2(t)\|^2-(r_{\mathrm{ob}}+r_2)^2, \; &t\in [0,\Tctrl]. 
\end{align*}

\begin{figure}[]
	\centering
	{\includegraphics[scale=0.95]{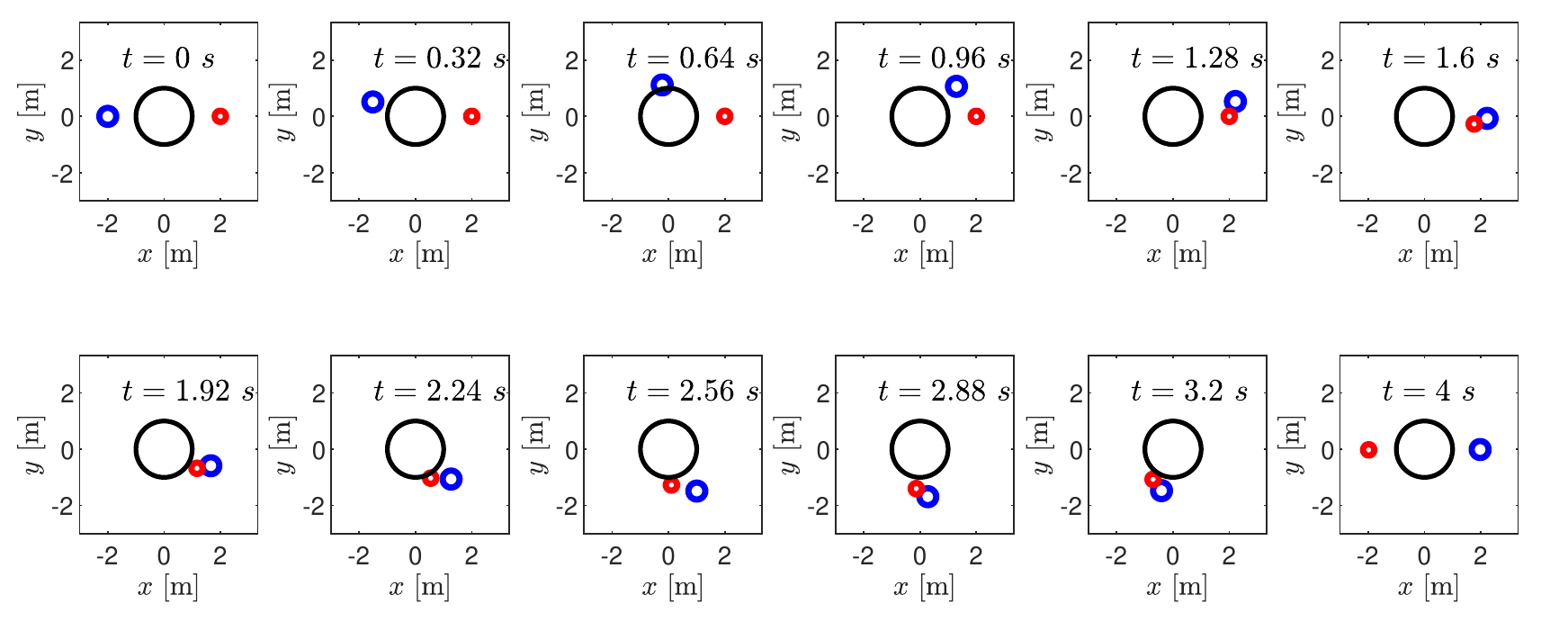}}
	\caption{
		Several frames of the optimal solution. 
		The first disc is marked with blue, the second with red, and the obstacle with black.
	}
	\label{fig:manipulation_frames}
\end{figure}
\begin{figure}[]
	\vspace{-0.2cm}
	\centering
	{\includegraphics[scale=0.95]{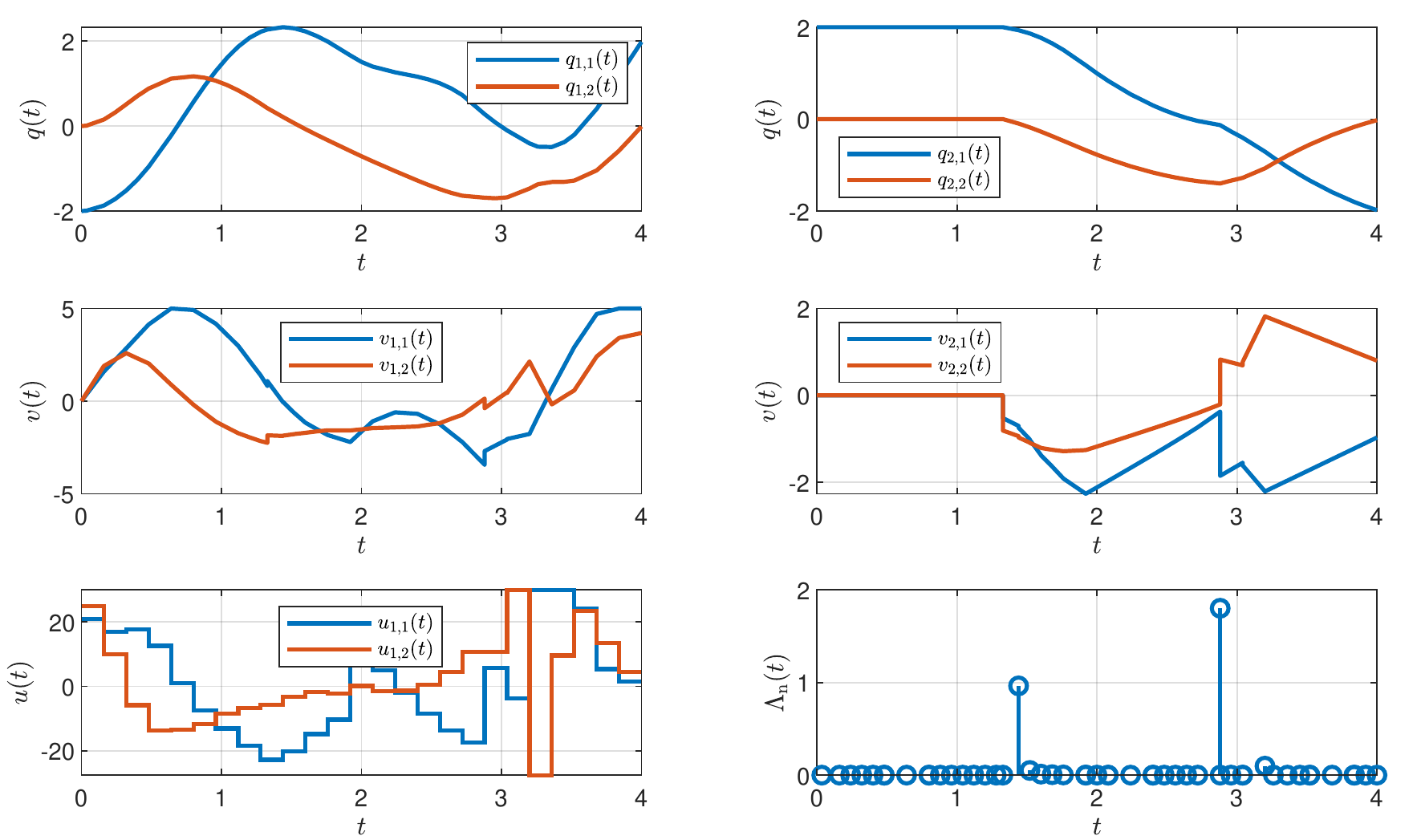}}
	\vspace{-0.2cm}
	\caption{The trajectories and optimal controls of the manipulation problem.}
	\label{fig:manipulation_states}
\end{figure}
\noindent The states are the positions of the two discs $q_1 = (q_{1,1},q_{1,2})$, $q_2 = (q_{2,1},q_{2,2})$ and their velocities $v_1 = (v_{1,1},v_{1,2})$, $v_2 = (v_{2,1},v_{2,2})$. 
The initial positions are $q_1(0) = (-2,0)$ and $q_2(0) = (2,0)$, and the initial velocities are zero $v_1(0)=v_2(0)=\mathbf{0}_{2,1}$.
The reference is $x_\mathrm{r} = (q_2(0),q_1(0),v_1(0),v_2(0))$, i.e., the discs should swap their positions and be at rest. 
Additionally, the overall control effort should be minimized.
This is modeled with least squares running and terminal objective terms with the weighting matrices:
$Q = \mathrm{diag}(10,10,10,10,0.01,0.01,0.01,0.01))$, $R = \mathrm{diag}(0.1,0.1)$ and $Q_T = 100Q$.
The control forces $u = (u_{1,1},u_{1,2})$ act only on the first disc.
We model a constant air friction force in the plane with $c_{\mathrm{f}}  = 2$ and $\epsilon_{\mathrm{f}} = 0.1$ to regularize the norm at zero.
The contact between the discs is frictionless.
The constant mass matrix reads as  $M = \mathrm{diag}(m_1,m_1,m_2,m_2)$ with $m_1 = 2$ and $m_2 = 1$.
The gap function reads as $f_c(q)=\|q_1(t)-q_2(t)\|^2-(r_1+r_2)^2$, where the radii of the balls are $r_1 = 0.3$ and $r_2 = 0.2$.
We introduce \textit{guiding} box constraints on the position and velocity and bound the control input in every direction.
The last two inequality constraints model the obstacle, which is a ball with a radius of $r_{\mathrm{ob}} = 1$ and located at the origin.

The problem is discretized and solved with NOSNOC. 
We discretize this OCP with a control horizon of $\Tctrl = 4$ with $\Nctrl = 25$ control intervals.
At every control interval, we the 3rd order Radau IIA \FESDj~method ($\Nstg = 2$) with $\NFE=2$.
The underlying MPCC is solved in a homotopy loop with a Scholtes relaxation \cite{Scholtes2001}, cf. \cite{Nurkanovic2022b} for more details.
The NLPs are solved with the interior-point solver IPOPT \cite{Waechter2006}, which is called via its CasADi interface \cite{Andersson2019}.

Figure \ref{fig:manipulation_frames} shows twelve frames of the optimal solution.
The optimizer finds a creative solution without any hints or sophisticated initial guesses.
The first ball goes to the second ball, makes contact with it, pushes it around the obstacle, and brings it to its final position.
It breaks the contact and returns to the final position.
Figure \ref{fig:manipulation_states} shows the states and optimal controls as a function of time.
One can see that the second ball is at rest until the first ball touches it and creates a velocity jump.

\FloatBarrier
\subsection{Optimal control problem with elastic and inelastic impacts}
We regard an optimal control problem with multiple contacts, of which some are elastic and some inelastic.
The problem involves three two-dimensional square boxes lying on a horizontal plane.
Only the middle box is actuated with a horizontal control force.
There are five possible contacts ($\ncontacts = 5$). 
The contact between the left and middle box is inelastic ($\CoR^1 =0 $) and between the middle and right is it partially elastic ($\CoR^2 = 0.5$).
Moreover, each box makes inelastic contact with the ground, i.e., $\CoRi = 0, i = 3,4,5$.
The middle box should, by making contact with the other two, bring all of them to their desired reference positions.
A one-dimensional variant without friction and similar to our example was considered in~\cite{Lin2022}.

We model this with the following optimal control problem:
\begin{subequations}\label{eq:ocp_box}
\begin{align}
	\min_{x(\cdot),\zn(\cdot),\zt(\cdot),u(\cdot),} \quad &  \int_{0}^{\Tctrl}\| x(t)-x_\mathrm{r}\|_Q^2+\|u(t)\|_R^2\;\dd t
	+\| x(\Tctrl)-x_\mathrm{r}\|_{Q_T}^2 \\
	\textrm{s.t.} \quad  &x(0) = {x}_0,\\
	&\textrm{Eq.}~\eqref{eq:cls}, \; t\in [0,\Tctrl], \\
	&q_\mathrm{lb} \leq q(t) \leq q_\mathrm{ub} , \; t\in [0,\Tctrl], \label{eq:ocp_box_guiding_q}\\
	&v_\mathrm{lb} \leq v(t) \leq v_\mathrm{ub} , \; t\in [0,\Tctrl], \label{eq:ocp_box_guiding_v}\\
	&u_\mathrm{lb} \leq u(t) \leq u_\mathrm{ub} , \; t\in [0,\Tctrl].
\end{align}
\end{subequations}
The states are the positions of the boxes: $q_1 = (q_{1,1},q_{1,2})$, $q_2 = (q_{2,1},q_{2,2})$, $q_3 = (q_{3,1},q_{3,2})$, collected into the vector $q=(q_1,q_2,q_3)\in\R^6$ and their velocities $v \in \R^6$, which are defined accordingly.
The initial positions are $q_1(0) = (-3,0)$, $q_2(0) = (0,0)$ and $q_3(0) = (3,0)$, which collected in the vector $q_0$. 
The initial velocities are zero $v_0=\mathbf{0}_{6,1}$, thus the initial state is $x_0 = (q_0,v_0)$.

All three boxes are square with side length $a = 2$, and they have the same masses of $m_i = 1, i =1,2,3$.
The inertia matrix is the constant diagonal matrix $M = \mathrm{diag}(m_1,m_1,m_2,m_2,,m_3,m_3)$.
The control forces $u(t) \in \R$ act only on the third box, i.e., $B_u =(0,0,1,0,0,0)$.
Beside the contact and friction forces, only gravity acts on the boxes, i.e.,
$f_v(q,v) = (0,-m_1 g,0,-m_2 g,0,-m_3 g)$ with $g =9.81$.
\begin{figure}
	\centering
	{\includegraphics[scale=0.85]{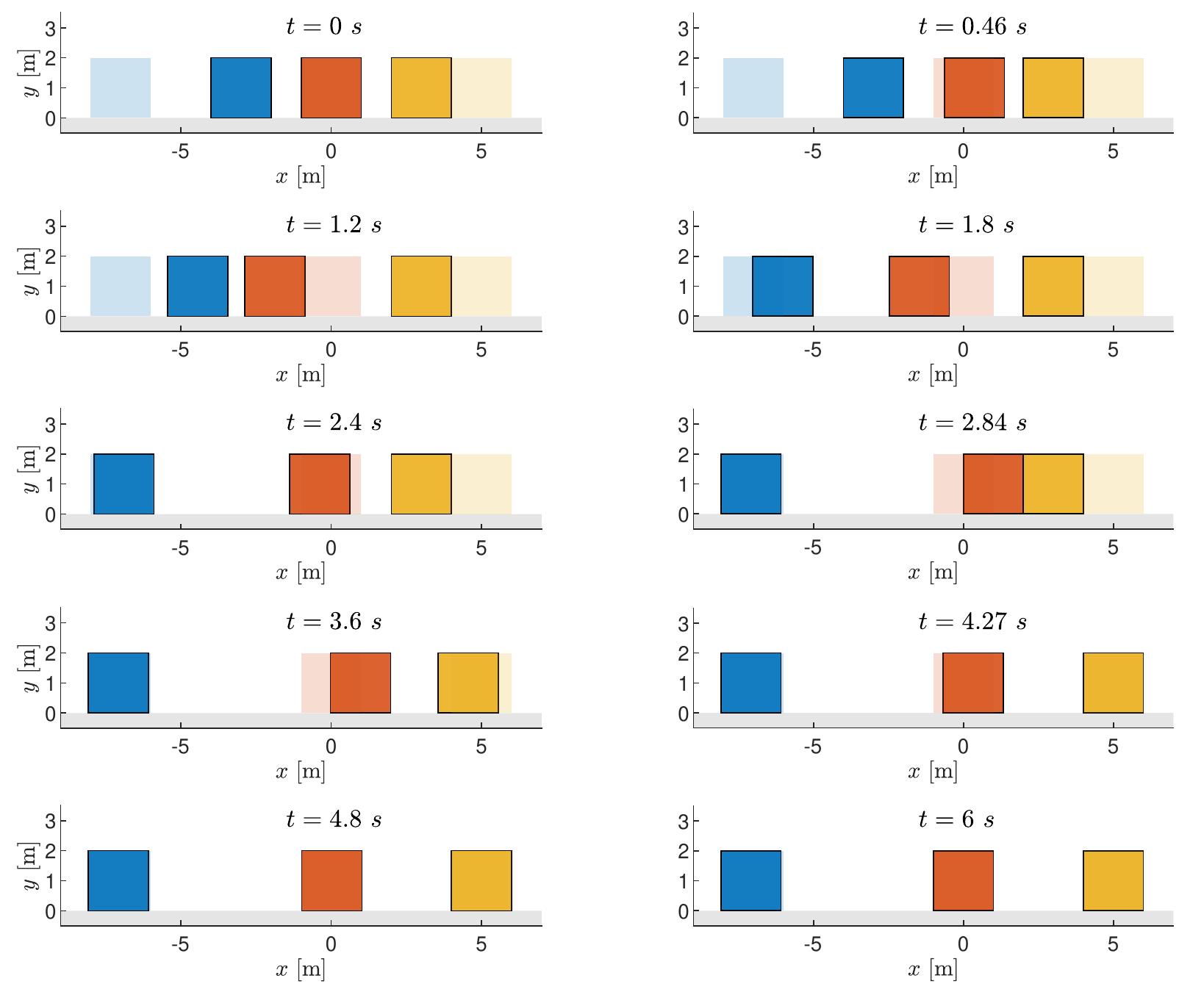}}
	\caption{
		Several frames of the optimal solution. 
		The shaded boxes show the desired final reference positions.
	}
	\label{fig:carts_frames}
\end{figure}
\begin{figure}
	\centering
	{\includegraphics[scale=0.8]{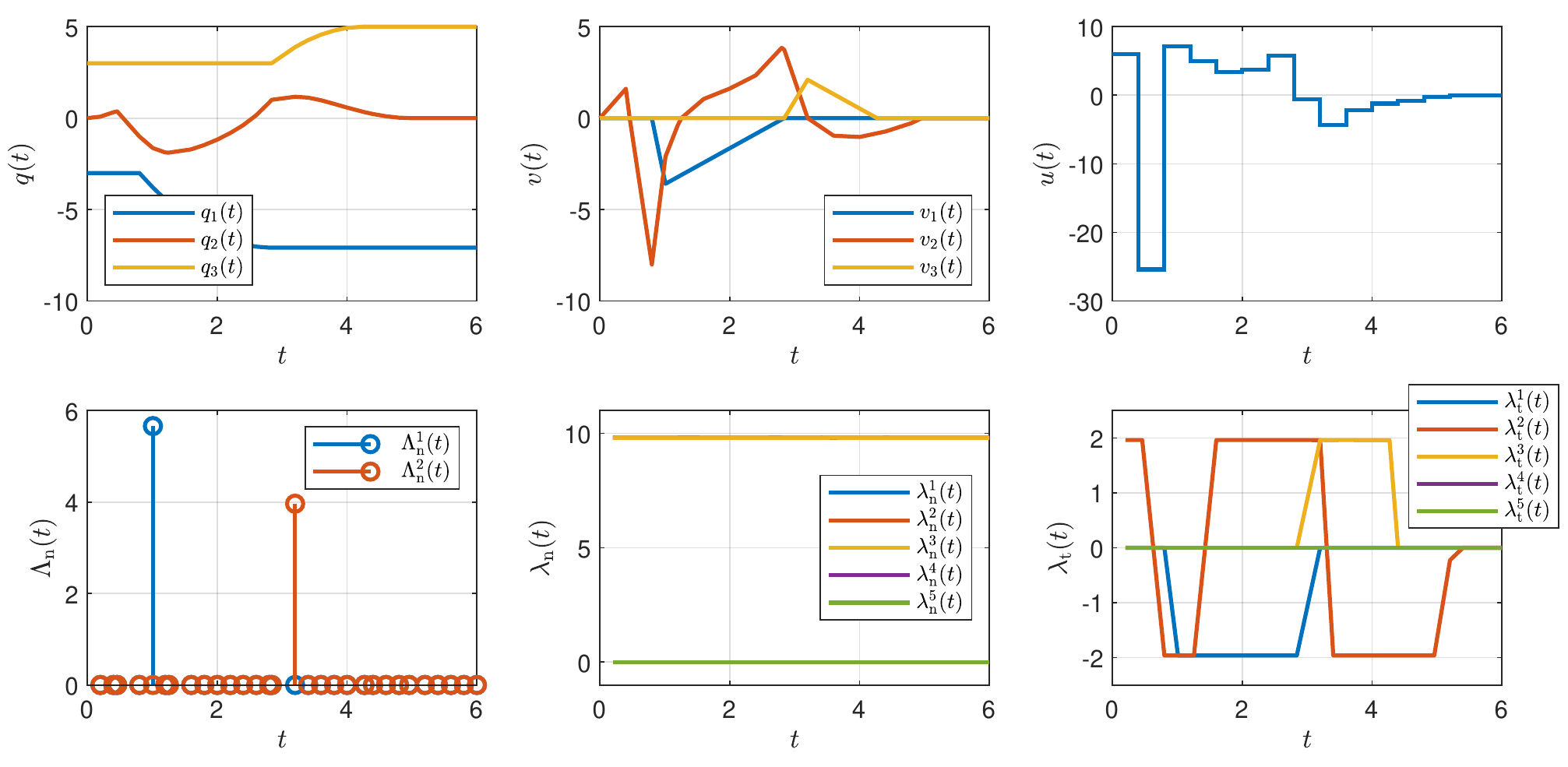}}
	\vspace{-0.2cm}
	\caption{The trajectories and optimal controls of the cart pushing problem.}
	\label{fig:carts_states}
\end{figure}
The gap function for the five possible contacts is:
\begin{align*}
	f_c(q) = \begin{bmatrix}
		q_{2,1} - q_{1,1} - a\\
		q_{3,1} - q_{2,1} - a\\
		q_{1,2} - \frac{a}{2}\\
		q_{2,2} - \frac{a}{2}\\
		q_{3,2} - \frac{a}{2}
	\end{bmatrix} \geq 0.
\end{align*}

All of the five contacts are frictional and the coefficient of friction are collected in the vector
$\mu = (0.1,0.1,0.2,0.2,0.2)$.
The control force is bounded and we have that $u_{\mathrm{lb}} = -30$ and $u_{\mathrm{ub}} = 30$.
Additionally, we add \textit{guiding constraints} \eqref{eq:ocp_box_guiding_q}-\eqref{eq:ocp_box_guiding_v}, which prevent and that the position and velocity take too large values in the intermediate iterates, which could impair the convergence.
We pick $ q_\mathrm{ub} = v_\mathrm{ub}  = 10e$ and $ q_\mathrm{lb} = v_\mathrm{lb}  = -10e$.

The desired final position are $q_\mathrm{r} = (-7,1,0,1,6,1)$, i.e., the left box should be moved further to the left, the right box further to right and the middle box should return to its initial position.
The boxes should rest at their final position, i.e., the reference velocity is $v_\mathrm{r} = \mathbf{0}_{6,1}$, and we define $x_{\mathrm{r}} = (q_{\mathrm{r}},v_{\mathrm{r}})$.
To achieve this motion, we define a least squares running and terminal objective terms with the weighting matrices:
$Q = \mathrm{diag}(10, 0.001, 1, 0.00, 10, 0.001,0.1, 0.1, 0.1, 0.1, 0.1, 0.1)$, $R = 0.1 $ and $Q_T = 100Q$.

The control horizon is $T = 6$, and the OCP \eqref{eq:ocp_box} is discretized using the \FESDj\ method with the 3rd order Radau IIA method ($\Nstg =2$), with $\Nctrl = 15$ equidistant control intervals and $\NFE =2$ finite elements on every control interval. 
The underlying MPCC is solved in a homotopy loop within NOSNOC~\cite{Nurkanovic2022b} using IPOPT~\cite{Waechter2006}.

Figure \ref{fig:carts_frames} shows ten frames of the optimal solution. 
The middle red box goes first slightly to the right, to have more space for accelerating. 
Afterwards, it hits the blue box which slides to its final position and stops due to friction. 
In the meantime, the middle red box hits the orange box on the right, which slides to its target position and stops due to the friction force.
The optimizer finds a creative solution without a sophisticated initial guess.
The solution contains two impacts and several transitions from stick to slip and slip to slick due to the friction forces.
Figure \ref{fig:carts_states} shows the states, optimal control, contact and friction forces and the impulses.
Note that $\lambda_{\mathrm{n}}^i = g$ for $ i = 3,4,5$ due to the contact with ground, and $\lambda_{\mathrm{n}}^i = 0$ for $i = 1,2$, since these contacts are only active for an atom of the time during the impacts.
The impact impulses are captured via $\Lambda_{\mathrm{n}}$.

\section{Discussion and outlook}\label{sec:conclusion}
This paper introduced the Finite Elements with Switch and Jump Detection (\FESDj) method for rigid bodies with friction and impact. 
This extends our previous work, where the FESD method was developed for ODEs with a discontinuous right-hand side~\cite{Nurkanovic2022b,Nurkanovic2023,Nurkanovic2022}.
\FESDj\ is an event-based discretization method that relies entirely on the solution of nonlinear complementarity problems and does not require any additional zero-location or mode selection algorithms.
Its four key algorithmic ingredients are: 
\begin{enumerate}
	\item The integration step sizes as degrees of freedom. 
	\item Cross complementarity conditions allow switches at the finite element boundary and implicitly force exact switch detection.
	\item Step equilibration conditions remove spurious degrees of freedom.
	\item Modified impact equations in complementarity form such that, depending on the active set, they impose either state jump or continuity conditions.
\end{enumerate}
The \FESDj~method can handle both elastic and inelastic impacts simultaneously. 
It can also handle frictional impacts and switch detection in slip-stick transitions.
It relies on the time-discretization of the nonlinear friction cone without polyhedral modeling approximations.
The method's intended use is for the direct transcription of continuous-time optimal control problems subject to a complementarity Lagrangian system (CLS).
We showed empirically that it recovers the integration accuracy that the underlying Runge-Kutta method has for smooth ODEs.
We recall that one can use any Runge-Kutta method within \FESDj.
Furthermore, we show that it is suitable for solving optimal control problems with both elastic and inelastic impacts, which provides great modeling flexibility.
Our work extends previous results \cite{Shield2022}, which was only used with Gauss-Legendre RK methods, used polyhedral friction cone approximations, was not covering all possible switching cases, and did not deal with spurious degrees of freedoms in the absence of switches.

The \FESDj\ method provides an alternative our previous approach, which transforms a CLS into a Filippov system via the time-freezing reformulation \cite{Nurkanovic2021a,Nurkanovic2021}, and uses the FESD method for discontinuous ODEs~\cite{Nurkanovic2022}.
Both approaches lead to high-accuracy approximations of the solution trajectory.
The advantage of FESD-J is that it can treat both elastic and inelastic impacts at the same time.
The discretization of an OCP with FESD or \FESDj~results in a Mathematical Program with Complementarity Constraints (MPCC).
We observed in our experiments that the MPCCs obtained from \FESDj~are difficult to solve with standard relaxation and penalty methods, i.e., where MPCCs are solved by solving a sequence of relaxed and more regular nonlinear programs with a standard solver.
In our experiments, we used IPOPT~\cite{Waechter2006} and SNOPT~\cite{Gill2005}, and observed similar performance.
In particular, for \FESDj we required more tuning of the homotopy parameters to obtain convergence for complicated OCPs compared to previous FESD approaches~\cite{Nurkanovic2021a,Nurkanovic2023}. 
Arguably, this is rather a limitation of the relaxation-based MPCC solutions strategies than the \FESDj\ formulation itself.
In future work, we aim to exploit pivoting-based MPCC methods~\cite{Leyffer2007,Kirches2022}, which explicitly treat the nonsmoothness and combinatorial structures in the complementarity constraints.
Moreover, with a more sophisticated MPCC solver at hand, we aim to perform an extensive numerical comparison of FESD with time-freezing and \FESDj.

Like any event-based method, FESD cannot, in principle, handle the accumulation of switching events, i.e., the Zeno phenomenon.
However, in an optimal control problem, more switches than finite elements would lead to an infeasible problem.
Therefore, if the optimizer finds an optimal solution, it will have finitely many switches.

It remains to rigorously prove the empirically observed high accuracy and convergence of numerical sensitivities.
We expect to obtain similar theoretical results as in \cite{Nurkanovic2022}, since FESD-J extends ideas of FESD for discontinuous ODEs.
In particular, we expect that the solutions to the FESD-J problems \eqref{eq:fesd_compact} have locally unique solutions, although they always consist of an overdetermined system of equations.
Without switches, the cross complementarities are implicitly satisfied by the standard stage-wise complementarity conditions. 
Similarly, if switches occur, the step equilibration conditions are implicitly satisfied.
Overall, it should be possible to formally show that for a fixed active set, the problem always reduces to a well-defined square system of equations.
Moreover, using the proof techniques from \cite{Nurkanovic2022}, one should be able to prove the convergence of solutions and numerical sensitivities of the \FESDj\ method with the same accuracy as the underlying RK method.

Interestingly, we observed in numerical experiments that if an accumulation of events is expected, e.g., in the bouncing ball example, the \FESDj\ method would result at some point in an impulse $\Zn$ that brings the ball to rest. 
It would be interesting to understand this event-capturing feature in more detail and to provide quantitative error estimates for this scenario.
Finally, another interesting extension would be to consider bilateral constraints within the complementarity Lagrangian model~\cite{Brogliato2016, Bruls2014}.

\bibliographystyle{plain}

\appendix

\section{Alternative formulation of impulse equations}\label{app:shield_equations}
	Inspired by \cite{Shield2022},  we rewrite the impulse equations in the complementarity form:
	\begin{subequations}\label{eq:impulse_eq_comp_shield}
		\begin{align}
			&M(q)(v(\tsp)-v(\tsn)) = \Jni(q)\Zni,\\
			&0\leq \Zni \perp \fci(q)  + |\Jni(q)^\top (v(\tsp) + \CoRi v(\tsn))| \geq 0 \label{eq:impulse_eq_comp_comp_shield}.
		\end{align}
	\end{subequations}
	Note that if $\fci(q) >0$ (no contact), we have that $\Zni =0$ and since the matrix $M(q)$ is full rank, we have that the velocity stays continuous, i.e., $v(\tsp)= v(\tsn)$. 
	We take the absolute value of the restitution law equations in \eqref{eq:impulse_eq_comp_comp_shield} since it can be negative when $\fci(q) >0$ and violate the nonnegativity condition.
	On the other hand, for  $\fci(q)=0$, we have that the state jumps are triggered, which implies that $\Zni >0$ and 
	$|\Jni(q(\ts))^\top (v(\tsp) + \CoRi v(\tsn))| = 0 \iff  \Jni(q(\ts))^\top (v(\tsp) + \CoRi v(\tsn))=0$, thus \eqref{eq:impulse_eq_plain} is recovered.
	The absolute value makes one of the functions involved in the complementarity constraint nondifferentiable.
	To get rid of this, we introduce another complementarity condition and rewrite \eqref{eq:impulse_eq_comp_shield} as
	\begin{subequations}\label{eq:impulse_eq_comp_smooth_shield}
		\begin{align}
			&M(q)(v(\tsp)-v(\tsn)) = \Jni(q)\Zni,\\
			&0\leq \Zni \perp \fci(q)  +\Slakcposi + \Slakcnegi \geq0,\\
			&\Slakcposi  - \Slakcnegi = \Jni(q)^\top (v(\tsp) + \CoRi v(\tsn)),\\
			&0\leq \Slakcposi  \perp \Slakcnegi \geq 0.
		\end{align}
	\end{subequations}
Except for the last complementarity condition, this formulation is equivalent to the one proposed in \cite{Shield2022}.
Note that without this complementarity condition, the variables $\Slakcposi,\Slakcnegi$ are not unique.
It is not difficult to see this formulation has similar properties at the one from Section~\ref{sec:impulse_equations}. 
However, it requires more variables	and mixes the complementarity conditions on position and velocity level.

\end{document}